\documentclass[preprints,article,accept,pdftex,moreauthors]{my-mdpi} 

\firstpage{1} 
\pubvolume{11}
\issuenum{7}
\articlenumber{1595}
\pubyear{2023}
\copyrightyear{2023}
\externaleditor{Academic Editors: António Lopes, Omid Nikan and Zakieh Avazzadeh}
\datereceived{17 January 2023} 
\daterevised{13 March 2023} 
\dateaccepted{23 March 2023} 
\datepublished{25 March 2023} 
\hreflink{https://doi.org/\\10.3390/math11071595}
\doinum{10.3390/math11071595}

\pdfoutput=1 

\Title{Three-Species Predator--Prey Stochastic Delayed Model Driven by L\'{e}vy Jumps and with Cooperation Among Prey Species}

\TitleCitation{Three-Species Predator–Prey Stochastic Delayed Model Driven by L\'{e}vy Jumps and with Cooperation Among Prey Species}

\Author{Jaouad Danane $^{1,\dagger}$\orcidA{} and 
Delfim F. M. Torres $^{2,}$*$^{,\dagger}$\orcidB{}}

\AuthorNames{Jaouad Danane and Delfim F. M. Torres}

\AuthorCitation{Danane, J.; Torres, D.F.M.}

\address{%
$^{1}$ \quad Laboratory of Systems, Modelization and Analysis for Decision Support,
National School of Applied Sciences (ENSA), University Hassan 1st of Settat, 
B.P 218 Berrechid, Morocco; jaouaddanane@gmail.com\\
$^{2}$ \quad Center for Research and Development in Mathematics and Applications (CIDMA),
Department of Mathematics, University of Aveiro, 3810-193 Aveiro, Portugal}

\corres{Correspondence: delfim@ua.pt; Tel.: +351-234-370-668}

\firstnote{These authors contributed equally to this work.} 

\abstract{{Our study focuses on analyzing the behavior of a stochastic 
predator–prey model with a time delay and logistic growth of prey, influenced 
by L\'{e}vy noise. Initially, we establish the existence, uniqueness, and boundedness 
of a positive solution that spans globally. Subsequently, we explore the conditions 
under which extinction occurs, and identify adequate criteria for persistence. Finally,
we validate our theoretical findings through numerical simulations, which also helps 
illustrate the dynamics of the stochastic delayed predator–prey model based on different criteria.}}

\keyword{prey–predator system; 
stochastic persistence and extinction;
logistic growth rate;
L\'{e}vy noise; 
time delays}

\MSC{34K50; 60H10}

\begin{document}

% ----------------------------------

\section{Introduction}

Recently, the prey–predator problem has attracted the interest of many researchers 
in ecology \cite{1,2,3}. The first prey–predator system, trying to describe 
the evolution of two sub-populations, the predator and the prey, was suggested 
in 1956 by Lotka and Volterra \cite{4}. This basic model has played a crucial 
starting view in different studies of prey–predator dynamics. Understanding 
the different reactions between prey–predator components becomes then an important 
issue to control the extinction of predator or prey. For example, several models 
are used to better understand the behavior of viral infections, e.g., in lions, 
wildebeest and zebra \cite{5}, and fishing \cite{6}. 

To better portray the growth of a population under constrained environments, 
Verhulst~\cite{7} introduced a special sort of growth equation, now known 
as the logistic equation. The logistic equation is very simple and very useful 
for one type of animal or one group of populations, as well as for
more than one population. Without a predator, the growth rate 
of the prey population fulfills a logistic growth equation. Much research has been 
done on modeling by distinguishing the parameters for various species 
\cite{8,9,10,11,12}. We should note that in some applications, such as the growth 
of population, the interaction of invulnerability, 
physiology, etc., means that the response growth rate of the species does not occur 
immediately, so it is necessary to include some kind of time delay \cite{MR4236351}. 
Delays have a significant effect on the stability of the system. In some examples, 
the system is stable for a small time-delay value and turns into an unstable system 
for a higher delay \cite{13}. The first to introduce a delay into the 
logistic differential equation was Hutchinson, with his paper from 1948 \cite{14}.

In the classical prey–predator model, it is assumed that the interaction 
is only within one prey and one predator. We can, however, take a model representing 
three interacting sub-populations, say one predator and two prey,
where the prey being modeled do not have competition among them, 
while mutual cooperation can be established between them 
against any predator. A work of 2018 by Kundu and Maitra proposes 
a prey–predator model with mutualism for two prey and one predator, also
taking into consideration the time delays, as follows \cite{15}:
\begin{equation}
\label{system1}
\begin{cases}
\dfrac{dx}{dt}(t)= r_1x(t)\left(1-\dfrac{x(t-\tau_1)}{K_1}\right)
-\alpha_1 x(t)z(t)+\beta x(t)y(t)z(t),\\
\dfrac{dy}{dt}(t)= r_2y(t)\left(1-\dfrac{y(t-\tau_2)}{K_2}\right)
-\alpha_2 y(t)z(t)+\beta x(t)y(t)z(t),\\
\dfrac{dz}{dt}(t) = -\delta z(t)-\alpha_3z^2(t)+a_1x(t-\tau_3)z(t)+a_2y(t-\tau_3)z(t),
\end{cases}
\end{equation}
where $x$ and $y$ are the densities of the prey, the predator is represented by $z$, 
and the time delays $\tau_1\geq 0$, $\tau_2\geq 0$ and $\tau_3\geq 0$ are
incorporated in the growth component for each species, 
subject to the initial data
\begin{adjustwidth}{-\extralength}{0cm}
\begin{equation}
\label{ID}
x(t)=\phi_1(t)\geq 0,\; 
y(t)=\phi_2(t)\geq 0, \; 
z(t)=\phi_3(t)\geq 0, \; 
t\in [-\tau,0],\; 
\tau=\max\{\tau_1,\tau_2, \tau_3\},
\end{equation}\vspace{-3pt}
\end{adjustwidth}
with $\phi_i\in C\left( [-\tau,0],\mathbb{R}_+\right)$ for $i=1,2,3$.
The description of the parameters of model \eqref{system1} 
is given in Table~\ref{tabl12}.
% ---------------------------------------
\begin{table}[H]
\caption{The description of the parameters for model \eqref{system1}.}
\label{tabl12} 
\newcolumntype{C}{>{\centering\arraybackslash}X}
\begin{tabularx}{\textwidth}{m{3cm}<{\centering}C}
\toprule
\textbf{Parameters} & \textbf{Description}\\
\midrule
$r_1$ 	& The intrinsic growth rate for prey $x$\\
$r_2$& The intrinsic growth rate for prey $y$\\
$K_1$& The carrying capacity for $x$\\
$K_2$& The carrying capacity for $y$\\
$\alpha_1$	& The predation rate of prey $x$ \\
$\alpha_2$	& The predation rate of prey $y$ \\		
$\beta$& The rate of cooperation of preys $x$ and $y$ against predator $z$\\
$\delta$& The predator death rate\\
$\alpha_3$	& The rate of intra-species competition within the predators\\
$a_1$ & The transformation rate of predator to preys $x$\\
$a_2$ & The transformation rate of predator to preys $y$\\
\bottomrule
\end{tabularx}
\end{table}
% ---------------------------------------

Many stochastic predator–prey problems have been deployed in order to describe 
the effect of white noise on predator–prey dynamics \cite{23,3,25,26,27}. 
In this context, in 2020, Rihan and Alsakaji studied
a stochastic prey–predator model by introducing a Brownian perturbation 
to the three variables of the model and obtaining different conditions 
of extinction and persistence \cite{23}. 

\textls[-15]{L\'evy jumps present an important tool to model many real dynamic 
phenomena~\cite{16,17,22,28,29,30,31}. Because of the inherent stochastic 
properties of the prey–predator processes, it is natural to assume that 
the dynamic model may experience sudden strong perturbations 
in the predator–prey progression \cite{32}. }

{The authors of the previous delay model in \cite{15} studied the effect 
of time delays on the stability of an equilibrium. Here we make model \eqref{system1} 
more realistic and improve the results of \cite{15} by considering that 
the predator and prey are exposed to uncertainties and randomness 
in the progress of a natural conflict. As a result, we consider here 
a stochastic system, which is more advantageous than the deterministic 
one because it takes into consideration not only the mean tendency 
but also the variance aspect surrounding it. In contrast to the deterministic 
model of \cite{15}, which always generates the same results for fixed initial values, 
our stochastic model may project different realities. In reality, here we add both 
random walks and L\'{e}vy jumps to \eqref{system1} and, in contrast to \cite{15}, 
we study the questions of persistence and extinction, while in \cite{15} their 
interests are different, being focused on bifurcation classification.} 
Motivated by these facts and by the previous studies, here, we introduce 
jumps into model \eqref{system1} as follows:
\begin{adjustwidth}{-\extralength}{0cm}
\begin{equation}
\label{sys2}
\begin{cases}
\displaystyle dx(t)
= \left(r_1x(t)\left(1-\dfrac{x(t-\tau_1)}{K_1}\right)
-\alpha_1 x(t)z(t)+\beta x(t)y(t)z(t)\right) dt +\sigma_1 x(t)dW_1(t)\\
\quad + \displaystyle\int_{U}q_1(u)x (t-)\tilde{N}(dt,du),\\
\displaystyle dy(t)= \left(  r_2y(t)\left(1-\dfrac{y(t-\tau_2)}{K_2}\right)
-\alpha_2 y(t)z(t)+\beta x(t)y(t)z(t)\right) dt +\sigma_2 y(t)dW_2(t)\\
\quad + \displaystyle\int_{U}q_2(u)y (t-)\tilde{N}(dt,du),\\
\displaystyle dz(t)
= \left(-\delta z(t)-\alpha_3z^2(t)+a_1x(t-\tau_3)z(t)
+a_2y(t-\tau_3)z(t)\right) dt +\sigma_3 z(t)dW_3(t)\\
\quad + \displaystyle\int_{U}q_3(u)z (t-)\tilde{N}(dt,du),
\end{cases}
\end{equation}
\end{adjustwidth}
where $W_i(t)$, $i = 1, 2, 3$, are standard Brownian motions defined on a complete 
probability space $\left( \Omega,\mathcal{F},(\mathcal{F}_t)_{t\geq 0},\mathbb{P}\right)$ 
with the filtration $(\mathcal{F}_t)_{t\geq 0}$ and satisfying the usual conditions. 
We denote by $x(t-)$, $y(t-)$  and $z(t-)$ the left limits of $x(t)$, $y(t)$ and $z(t)$, 
respectively. Here $N(dt,du)$ is a Poisson counting measure with the stationary 
compensator $\nu(du)dt$, $\tilde{N}(dt,du)=N(dt,du)-\nu(du)dt$, where $\nu$ is 
defined on a measurable subset $U$ of the nonnegative half-line with 
$\nu(U)< \infty$ and the intensity of $W_i(t)$ is $\sigma_i$, $i=1,2,3$. 
The jump intensities are represented by $q_i(u)$, $i=1,2,3$.

Our model and work is new and original
because it introduces both white and L\'{e}vy noise into the system \eqref{system1}, 
to show the effect on the sub-population of the prey 
under their cooperation and by taking their recruitment rate in logistic growth,   
and considers the effect of the stochastic perturbation of the predator. On the other hand, 
we also deal for the first time in the literature with the effect of the delays 
on the three considered species, describing the intra-production of the new prey
and the generation of the predator after reaction with the prey.

The paper is organized as follows. In Section~\ref{sec:2}, 
we give some properties of the solution of model \eqref{sys2}. 
The stochastic extinction of both the prey and the predator is established in Section~\ref{sec:3}. 
In Section~\ref{sec:4}, we prove the persistence of the prey and the extinction of the predator. 
The stochastic persistence of both the predator and the prey is obtained 
in Section~\ref{sec:5} and Section~\ref{sec:6} is devoted to some numerical simulations 
that illustrate the theoretical findings. We end with Section~\ref{sec:7} 
of conclusions, final discussions, and some possible 
future directions of research.

% ----------------------------------

\section{Properties of the Solution}
\label{sec:2}

In this section, we show the well-posedness 
of our stochastic predator–prey delayed model
driven by L\'{e}vy jumps. In reality, we prove the 
existence and uniqueness of a global positive solution 
(Theorem~\ref{thm1}) as well as its boundedness 
(Theorem~\ref{thm:AS1}).

% ----------------------------------

\subsection{Existence and Uniqueness of a Global Positive Solution}

Our first theorem guarantees the existence and uniqueness 
of a global positive solution for \eqref{sys2}.
Please note that the unique global solution 
$\left(x(t),y(t),z(t)\right)$ predicted
by Theorem~\ref{thm1} satisfies 
$x(t)\leq K_1$ and $y(t)\leq K_2$
since $K_1$ is the carrying capacity 
for $x$ and $K_2$ is the carrying capacity for $y$.

\begin{Theorem}
\label{thm1}
Let $\delta>\alpha_3$. For any given initial condition \eqref{ID},
the system \eqref{sys2} has a unique global solution 
$\left(x(t),y(t),z(t)\right) \in \mathbb{R}_+^3$,
$t\geq -\tau$, almost surely (a.s.).
\end{Theorem}

\begin{proof}
The drift and the diffusion are locally Lipschitz. Thus, 
for any initial condition~\eqref{ID}, there exists a unique local solution  
$\left(x(t),y(t),z(t)\right)$ for $t\in [-\tau,\tau_e)$, where $\tau_e$ is the explosion time.
To demonstrate that this solution is global, we need to prove that $\tau_e=\infty$ a.s. 
First, we will prove that $\left(x(t),y(t),z(t)\right)$ do not tend to infinity for a finite time. 
Let $k_0>0$ be a sufficiently large number, in such a way that 
$\left(\phi_1,\phi_2,\phi_3\right)\in C\left([-\tau,0],\mathbb{R}^3_+\right)$ 
is within the interval $[1/k_0,k_0]$. Let us define, for each integer $k\geq k_0$, the stopping time
$$
\tau_k=\inf\left\{t\in[0,\tau_e) : x(t)\notin(1/k,k) 
\text{ or } 
y(t)\notin(1/k,k) 
\text{ or }  
z(t)\notin(1/k,k)\right\}.
$$
We have established that $\tau_k$ is an increasing number when $k\uparrow \infty$. 
Let $\tau_\infty=\lim_{k\rightarrow\infty}\tau_k$, where
$\tau_\infty\leq \tau_e$ a.s. We need to show that $\tau_\infty=\infty$, 
which means that $\tau_e=\infty$ and $\left(x(t),y(t),z(t)\right)\in \mathbb{R}_+^3$ a.s. 
Assume the opposite is verified, i.e., assume that $\tau_\infty<\infty$ a.s. 
Then, there exist two constants $0<\epsilon<1$ and $T>0$ such that
$\mathbb{P}(\tau_\infty\leq T)\geq \epsilon$. Consider the functional
\begin{align*}
V\left(x(t),y(t),z(t)\right)
=&\left( x-1-\ln\left(x\right)\right)+\dfrac{r_1}{K_1}\int_t^{t+\tau_1}x(s-\tau_1)ds\\
&+\left( y-1-\ln\left(y\right)\right)+\dfrac{r_2}{K_2}\int_t^{t+\tau_1}y(s-\tau_1)ds\\
&+\left( z-1-\ln\left(z\right)\right)+\dfrac{\rho}{2}
\int_t^{t+\tau_1}\left( a_1x^2(s-\tau_3)+a_2y^2(s-\tau_3)\right) ds,
\end{align*}
where $\rho=\dfrac{2\beta K_1K_2+\alpha_1+\alpha_2+a_1K_1+a_2K_2}{\delta-\alpha_3}$. 
From It\^{o}'s formula \cite{322}, we have
\begin{adjustwidth}{-\extralength}{0cm}
\begin{equation}
\label{GS}
\begin{split}
&dV\left(x,y,z\right)
= LV\;dt+\sigma_1(x-1)\;dW1+\sigma_{2}(y-1)\;dW_2+\sigma_3(z-1)\;dW_{3}\\
&+\int_U\left[ q_1(u)x-\log\left(1+q_1(u)x\right)+\dfrac{r_1}{K_1}q_1(u)
\int_t^{t+\tau_1}x(s-\tau_1)ds\right]\tilde{N}(dt,du)\\
&+\int_U\left[ q_2(u)y-\log\left(1+q_2(u)y\right)+\dfrac{r_2}{K_2}q_2(u)
\int_t^{t+\tau_2}y(s-\tau_2)ds\right]\tilde{N}(dt,du)\\
&+\int_U\left[q_3(u)z-\log\left(1+q_3(u)z\right)+\dfrac{a_1\rho}{2}
\int_t^{t+\tau_3}\left(2q_1(u)x(s-\tau_3)+q_1^2(u)x^2(s-\tau_3)\right) ds\right]\tilde{N}(dt,du)\\
&+\int_U\left[ \dfrac{a_2\rho}{2}\int_t^{t+\tau_3}\left(q_2(u)y(s-\tau_3)
+q_2^2(u)y^2(s-\tau_3)\right) ds\right]\tilde{N}(dt,du),
\end{split}
\end{equation}
\end{adjustwidth}
where
\begin{align*}
LV=&\left(1-\dfrac{1}{x}\right)\left( r_1x(t)\left(1-\dfrac{x(t-\tau_1)}{K_1}\right)-\alpha_1x(t)z(t)
+\beta x(t)y(t)z(t)\right)+\dfrac{\sigma_1^2}{2}\\
&+\left(1-\dfrac{1}{y}\right)\left( r_2y(t)\left(1-\dfrac{y(t-\tau_2)}{K_2}\right)-\alpha_2y(t)z(t)
+\beta x(t)y(t)z(t)\right)+\dfrac{\sigma_2^2}{2}\\
&+\dfrac{\rho a_1}{2}\left(x^2-x^2(t-\tau_3)\right)+\dfrac{\rho a_2}{2}\left(y^2-y^2(t-\tau_3)\right)\\
&+\rho\left(1-\dfrac{1}{z}\right)\left(-\delta z(t)-\alpha_3z^2(t)
+a_1x(t-\tau_3)z(t)+a_2y(t-\tau_3)z(t)\right)+\dfrac{\rho\sigma_3^2}{2}\\
&+\int_U\left[ q_1(u)-\log\left(1+q_1(u)\right)\right]\nu(du)\\
&+\int_U\left[ q_2(u)-\log\left(1+q_2(u)\right)\right]\nu(du)\\
&+\int_U\left[ q_3(u)-\log\left(1+q_3(u)\right)\right]\nu(du)\\
&+\dfrac{\rho}{2}\int_U\left[ q^2_1(u)\int_t^{t+\tau_3}a_1x^2(s-\tau_3)ds
+q^2_2(u)\int_t^{t+\tau_3}a_2y^2(s-\tau_3)ds\right]\nu(du).
\end{align*}
Therefore, we have
\begin{align*}
LV \leq& r_1x+2\beta xyz+r_1\dfrac{x(t-\tau_1)}{K_1}+\alpha_1z+\dfrac{\sigma_1^2}{2}\\
&+r_2y+r_2\dfrac{y(t-\tau_2)}{K_2}+\alpha_2z+\dfrac{\sigma_2^2}{2}
+\dfrac{\rho a_1}{2}x^2(t-\tau_3)+\dfrac{\rho a_2}{2}y^2(t-\tau_3)\\
&+\rho \delta+\rho \alpha_3 z-\rho \delta z+a_1xz+a_2yz+4C',\\
\end{align*}
\vspace{-24pt}
\begin{align*}
LV\leq& r_1K_1+2\beta K_1K_2z+r_1+\alpha_1z+\dfrac{\sigma_1^2}{2}+\dfrac{\sigma_3^2}{2}\\
&+r_2K_2+r_2+\alpha_2z+\dfrac{\sigma_2^2}{2}+\dfrac{\rho a_1}{2}K_1^2+\dfrac{\rho a_2}{2}K_2^2\\
&+\rho \delta+\rho \alpha_3 z-\rho \delta z+a_1K_1z+a_2K_2z+4C',\\
LV\leq & r_1K_1+\rho \delta+r_1+\dfrac{\sigma_1^2}{2}+\dfrac{\sigma_3^2}{2}\\
&+r_2K_2+r_2+\dfrac{\sigma_2^2}{2}+\dfrac{\rho a_1}{2}K_1^2+\dfrac{\rho a_2}{2}K_2^2\\
&+\underbrace{\left(2\beta K_1K_2+\alpha_1+\alpha_2+a_1K_1+a_2K_2
+\rho \alpha_3 -\rho \delta\right)}_{=0} z+4C',\\
LV\leq & C,
\end{align*}
where
\begin{align*}
C'=&\max\bigg\{\int_U\left( q_1(u)-\log\left(1+q_1(u)\right)\right)\nu(du),
\int_U\left( q_2(u)-\log\left(1+q_2(u)\right)\right)\nu(du),\\
&\int_U\left(q_3(u)-\log\left(1+q_3(u)\right)\right)\nu(du),
\dfrac{\rho \tau_3}{2} \int_U\left(q^2_1(u)a_1K_1^2+q^2_2(u)a_2K_2^2\right)\nu(du)\bigg\},
\end{align*}
and 
\begin{equation*}
C= r_1K_1+\rho \delta+r_1+\dfrac{\sigma_1^2}{2}+\dfrac{\sigma_3^2}{2}
+r_2K_2+r_2+\dfrac{\sigma_2^2}{2}+\dfrac{\rho a_1}{2}K_1^2+\dfrac{\rho a_2}{2}K_2^2+4C'.
\end{equation*}
Integrating both sides of equation \eqref{GS} between $0$ and $\tau_k\wedge T$, we obtain
\begin{equation*}
\begin{split}
0&\leq \mathbb{E}\left(V\left(x(\tau_k\wedge T),y(\tau_k\wedge T),z(\tau_k\wedge T)\right)\right)\\
&\leq V\left(x(0),y(0),z(0)\right)+CT.
\end{split}
\end{equation*}
For each $h>0$, let us define
$$
H(h)=\inf\left\{V(x_1,x_2,x_3) : x_i\geq h \text{ or } x_i\leq \dfrac{1}{h},\; i=1,2,3\right\},
$$
where $x_1=x$, $x_2=y$ and $x_3=z$. Then, 
$$
\lim_{h\rightarrow\infty}H(h)=\infty.
$$
Consequently, letting $k\rightarrow\infty$, we have
$$
\infty > V\left(x(0),y(0),z(0)\right)+CT =\infty, 
$$
which represents a contradiction of the previous assumption. 
Therefore, $\tau_\infty = \infty$ and the model has a unique global 
solution $(x(t),y(t),z(t))$ a.s.
\end{proof}

% ----------------------------------

\subsection{Stochastic Boundedness}

Theorem~\ref{thm1} shows that the solution of system \eqref{sys2} 
remains in the positive cone $\mathbb{R}_+^3$. However, this nonexplosion 
property in a population dynamical system is often not good enough: 
the property of ultimate boundedness is more desired. Next, we prove
stochastically ultimate boundedness.

\begin{Theorem}
\label{thm:AS1}
Under the conditions
\begin{equation}
\label{AS1}
\begin{gathered}
\sigma_1^2+\int_Uq_1^2(u)\nu(du)+2r_1+\beta K_2-\alpha_1 K_1<0,\\
\sigma_2^2+\int_Uq_2^2(u)\nu(du)+2r_2+\beta K_1-\alpha_2 K_2<0,\\
\sigma_3^2+\int_Uq_3^2(u)\nu(du)+2a_1K_1+2a_2K_2-\delta-\alpha_1 K_1-\alpha_2 K_2<0,
\end{gathered}
\end{equation}
the solution of system \eqref{sys2} given by Theorem~\ref{thm1}
is stochastically ultimately bounded for any initial condition \eqref{ID}.
\end{Theorem}

\begin{proof}
Consider function
$V\left(x,y,z\right)=x^2+y^2+z^2$. Using It\^{o}'s formula, we have 
\begin{equation}
\label{SBD}
\begin{split}
dV\left(x,y,z\right)
=&LV\; dt+2\sigma_1^2x^2dW_1(t)+2\sigma_2^2y^2dW_2(t)+2\sigma_3^2z^2dW_3(t)\\
&+\int_U(2q_1(u)+q_1^2(u))x^2\tilde{N}(dt,du)+\int_U(2q_2(u)+q_2^2(u))y^2\tilde{N}(dt,du)\\
&+\int_U(2q_3(u)+q_3^2(u))z^2\tilde{N}(dt,du),
\end{split}
\end{equation}
where
\begin{align*}
LV=
&2x\left(r_1x\left(1-\dfrac{x(t-\tau_1)}{K_1}\right)-\alpha_1xz
+\beta xyz\right)+2y\left(r_2y\left(1-\dfrac{y(t-\tau_2)}{K_2}\right)
-\alpha_2yz+\beta xyz\right)\\
&+2z\left(-\delta z-\alpha_3 z^2+a_1x(t-\tau_3)z+a_2y(t-\tau_3)z\right)\\
&+\sigma_1x^2+\sigma_2y^2+\sigma_3z^2+\int_Uq_1^2x^2\nu(du)\\
&+\int_Uq_2^2y^2\nu(du)+\int_Uq_3^2z^2\nu(du).
\end{align*}
Then,
\begin{align*}
LV
\leq& \left(\sigma_1^2+\int_Uq_1^2(u)\nu(du)+2r_1+\beta K_2-\alpha_1K_1\right)x^2\\
&+\left(\sigma_2^2+\int_Uq_2^2(u)\nu(du)+2r_2+\beta K_1-\alpha_2K_2\right)y^2\\
&+\left(\sigma_3^2+\int_Uq_3^2(u)\nu(du)+2a_1K_1+2a_2K_2-\delta-\alpha_1K_1-\alpha_2K_2\right)z^2.
\end{align*}
Put
\begin{align*}
f(x,y,z) 
=&\left(\sigma_1^2+\int_Uq_1^2(u)\nu(du)+2r_1+\beta K_2-\alpha_1K_1\right)x^2\\
&+\left(\sigma_2^2+\int_Uq_2^2(u)\nu(du)+2r_2+\beta K_1-\alpha_2K_2\right)y^2\\
&+\left(\sigma_3^2+\int_Uq_3^2(u)\nu(du)+2a_1K_1+2a_2K_2-\delta-\alpha_1K_1-\alpha_2K_2\right)z^2.
\end{align*}
Using conditions \eqref{AS1}, we find that function $f(x,y,z)$ has an upper bound. Denote 
$$
M=\sup_{(x,y,z)\in \mathbb{R}^3_+}f(x,y,z) 
\text{ and } N=M+1.
$$
Since $f(0,0,0) = 0$, then $N>0$. According to formula \eqref{SBD}, 
\begin{align*}
dV 
\leq& \left[N-(x^2+y^2+z^2)\right]dt+2\sigma_1^2x^2dW_1(t)
+2\sigma_2^2y^2dW_2(t)+2\sigma_3^2z^2dW_3(t)\\
&+\int_U(2q_1(u)+q_1^2(u))x^2\tilde{N}(dt,du)
+\int_U(2q_2(u)+q_2^2(u))y^2\tilde{N}(dt,du)\\
&+\int_U(2q_3(u)+q_3^2(u))z^2\tilde{N}(dt,du).
\end{align*}
Then, using It\^{o}'s formula, we obtain that
\begin{equation}
\label{INT1}
\begin{split}
d\left[e^tV\right] =&e^tV\; dt+e^tdV\\
\leq & Ne^t\;dt+2\sigma_1^2x^2dW_1(t)
+2\sigma_2^2y^2dW_2(t)+2\sigma_3^2z^2dW_3(t)\\
&+\int_U(2q_1(u)+q_1^2(u))x^2\tilde{N}(dt,du)
+\int_U(2q_2(u)+q_2^2(u))y^2\tilde{N}(dt,du)\\
&+\int_U(2q_3(u)+q_3^2(u))z^2\tilde{N}(dt,du).
\end{split}
\end{equation}
Integrating both sides of \eqref{INT1} from $0$ to $t$, 
and then taking expectations, we have
\begin{equation*}
e^t\mathbb{E}[V(X)] 
\leq V\left(x(0),y(0),z(0)\right)+ Ne^t-N
\qquad \forall t\in [-\tau,0],
\end{equation*}
where $X=(x,y,z)$. This fact implies that 
\begin{equation*}
\limsup_{t\rightarrow \infty}\mathbb{E}[V(X)] \leq N.
\end{equation*}
Since $V(X)=x^2+y^2+z^2$, then 
$$
\limsup_{t\rightarrow \infty}\mathbb{E}[|(X)|^2] \leq N.
$$
For any $\epsilon>0$, let $A=\dfrac{\sqrt{N}}{\sqrt{\epsilon}}$.  
Using Chebyshev's inequality, we obtain that
$$
\mathbb{P}\left( |(X)|>A\right)\leq \dfrac{\mathbb{E}[|(X)|^2] }{A^2}
\leq \dfrac{N}{\dfrac{N}{\epsilon}}=\epsilon.
$$
The proof is complete.
\end{proof}

We remark that our condition \eqref{AS1} assures the boundedness of the solution, 
which is crucial in biological models. Roughly speaking, Theorem~\ref{thm:AS1}
tells us that we need to control the intensities $\sigma_i$ of the Brownian motions
and the jump intensities $q_i$.

% ----------------------------------

\section{Stochastic Extinction of the Prey and Predator}
\label{sec:3}

In this section, we show conditions under which 
the population becomes extinct with probability of one
(Theorem~\ref{thm:extinction:prey:predator}).
In the follow-up, we use the following notation:
$$
\langle x(t)\rangle := \dfrac{1}{t}\int_0^t x(s)ds.
$$

\begin{Theorem}
\label{thm:extinction:prey:predator}
Let
\begin{align*}
c_1&=r_1-\dfrac{\sigma_1^2}{2},\\
c_2&=r_2-\dfrac{\sigma_2^2}{2},\\
c_3&=a_1\dfrac{K_1}{r_1}\left( r_1-\dfrac{\sigma_1^2}{2}\right)
+a_2\dfrac{K_2}{r_2}\left( 
r_2-\dfrac{\sigma_2^2}{2}\right)-\delta -\dfrac{\sigma_3^2}{2}.
\end{align*} 
If $\max\{c_1,c_2,c_3\}<0$, then the solution of system \eqref{sys2} satisfies
\begin{align*}
&\limsup_{t\rightarrow +\infty} \left<x(t)\right> =0,\\
&\limsup_{t\rightarrow +\infty} \left<y(t)\right> =0,\\
&\limsup_{t\rightarrow +\infty} \left<z(t)\right> =0,
\end{align*}
for any initial condition \eqref{ID}.
\end{Theorem}

\begin{proof}
Let us define
$$ 
F\left(x\right)=\log\left(x(t)\right)
-\dfrac{r_1}{K_1}\int_t^{t+\tau_1}x(s-\tau_1)ds.
$$
Using It\^{o}'s formula, we have 
\begin{equation*}
\begin{aligned}
dF=&LF\; dt+\sigma_1dW_1(t)
+\int_U \log\left(1+q_1(u)\right)\tilde{N}(dt,du),
\end{aligned}
\end{equation*}
where
\begin{equation*}
LF=r_1-\dfrac{\sigma_1^2}{2}+\int_U\log\left(1+q_1(u)\right)-q_1(u)\nu(du)
-\dfrac{r_1x(t)}{K_1}-z(t)\left(\alpha_1-\beta y(t)\right).
\end{equation*}
Observe that $x \mapsto \log(1+x)-x$ is a nonpositive function. Since
$$
LF\leq  r_1-\dfrac{\sigma_1^2}{2}-\dfrac{r_1x(t)}{K_1},\\
$$
then
\begin{align*}
&\dfrac{1}{t}\left[\log\left(x(t)\right)-\dfrac{r_1}{K_1}
\displaystyle \int_t^{t+\tau_1}x(s-\tau_1)ds-\log\left(x(0)\right)
+\dfrac{r_1}{K_1}\int_0^{\tau_1}x(s-\tau_1)ds\right]\\
&\leq r_1-\dfrac{\sigma_1^2}{2}-\dfrac{r_1}{K_1}\left<x(t)\right>
+\int_0^t\dfrac{\sigma_1dW_1(s)}{t}ds+\dfrac{1}{t}
\int_0^t \int_U \log\left(1+q_1(u)\right)\tilde{N}(ds,du)ds,
\end{align*}
and so
\begin{align*}
\left<x(t)\right>
\leq &  \dfrac{K_1}{r_1}\left( r_1-\dfrac{\sigma_1^2}{2}\right)
+\dfrac{K_1}{r_1} \left(\int_0^t\dfrac{\sigma_1dW_1(s)}{t}ds
+\dfrac{1}{t}\int_0^t \int_U \log\left(1+q_1(u)\right)\tilde{N}(ds,du)ds\right)\\
&-\dfrac{K_1}{r_1}\left(\dfrac{\log\left(x(t)\right)-\dfrac{r_1}{K_1}
\int_t^{t+\tau_1}x(s-\tau_1)ds-\log\left(x(0)\right)
+\dfrac{r_1}{K_1}\int_0^{\tau_1}x(s-\tau_1)ds}{t}\right).
\end{align*}
Put
$$
M_t=\int_0^t \sigma_1 dW_1(s).
$$
Thus,
\begin{align*}
\limsup_{t\rightarrow+\infty} \dfrac{\left<M_t,M_t\right>}{t}
=\limsup_{t\rightarrow+\infty} \dfrac{\sigma_1^2}{t}<\infty.
\end{align*}
Now, using the strong law of large numbers for martingales, we obtain
\begin{align*}
\liminf_{t\rightarrow+\infty}\dfrac{M_t}{t}=0.
\end{align*}
Moreover, one has
$$
\int_t^{t+\tau_1}x(s-\tau_1)ds
=\int_{t-\tau_1}^tx(s)ds
=\int_0^tx(s)ds-\int_0^{t-\tau_1}x(s)ds,
$$
and therefore
$$
\lim_{t\rightarrow +\infty}\dfrac{1}{t}\int_t^{t+\tau_1}x(s-\tau_1)ds
=\lim_{t\rightarrow +\infty}\dfrac{1}{t}
\left(\int_0^tx(s)ds-\int_0^{t-\tau_1}x(s)ds\right)=0.
$$
In addition, 
$$
\lim_{t\rightarrow +\infty}\dfrac{1}{t}\int_{-\tau_1}^0\phi_1(s)ds=0,
$$
so
\begin{equation*}
\limsup_{t\rightarrow +\infty}\left<x(t)\right>
\leq \dfrac{K_1}{r_1}\left( r_1-\dfrac{\sigma_1^2}{2}\right).
\end{equation*}
Using the same method, we prove that 
\begin{equation*}
\limsup_{t\rightarrow +\infty}\left<y(t)\right>
\leq \dfrac{K_2}{r_2}\left( r_2-\dfrac{\sigma_2^2}{2}\right).
\end{equation*}
Now, consider 
$$
V(t)=\log(z)+\int_t^{t+\tau_3}\left(a_1x(s-\tau_3)+a_2y(s-\tau_3)\right) ds.
$$
Using It\^{o}'s formula, it follows that
\begin{equation*}
dV=LV\; dt+\sigma_3dW_3(t)+\int_U \log\left(1+q_3(u)\right)\tilde{N}(dt,du),
\end{equation*}
where
\begin{align*}
LV=&-\delta-\alpha_3z+a_1x+a_2y-\dfrac{\sigma_3^2}{2}
+\int_U\log\left(1+q_3(u)\right)-q_3(u)\nu(du)\\
\leq &-\delta-\alpha_3z+a_1x+a_2y-\dfrac{\sigma_3^2}{2}.
\end{align*}
Then,
\begin{align*}
\dfrac{V(t)-V(0)}{t}
\leq & a_1\left<x\right>+a_2\left<y\right>-\delta 
-\alpha_3\left<z\right>-\dfrac{\sigma_3^2}{2}\\
&+\dfrac{1}{t}\int_0^t \left(\sigma_3dW_3(t)
+\int_U \log\left(1+q_3(u)\right)\tilde{N}(dt,du)\right),
\end{align*}
and
\begin{align*}
\left<z\right>
\leq &\dfrac{1}{\alpha_3}\left(a_1\left<x\right>+a_2\left<y\right>
-\delta -\dfrac{\sigma_3^2}{2}-\dfrac{V(t)-V(0)}{t}\right)\\
&+\dfrac{1}{\alpha_3}\dfrac{1}{t}\int_0^t \left(\sigma_3dW_3(t)
+\int_U \log\left(1+q_3(u)\right)\tilde{N}(dt,du)\right).
\end{align*}
Using the same technique we have used for $\left<x\right>$, we obtain that
$$
\limsup_{t\rightarrow +\infty}\left<z\right>  
\leq \dfrac{1}{\alpha_3}\left(a_1\dfrac{K_1}{r_1}\left( r_1
-\dfrac{\sigma_1^2}{2}\right)+a_2\dfrac{K_2}{r_2}\left( 
r_2-\dfrac{\sigma_2^2}{2}\right)-\delta -\dfrac{\sigma_3^2}{2}\right).
$$
This completes the proof.
\end{proof}

% ----------------------------------

\section{Stochastic Extinction of Predator}
\label{sec:4}

We now show conditions for which the population 
of predators becomes extinct with probability of one
(Theorem~\ref{thm:exct:persist}).

\begin{Definition}
A positive function $\mathcal{V}$ 
is said to be persistent in mean if 
\begin{equation*}
\liminf_{t \to \infty}\inf\left<\mathcal{V}\right>
=\liminf_{t \to \infty}\inf\frac{1}{t}\int_{0}^{t} \mathcal{V}(u)du >0.
\end{equation*}
\end{Definition}

\begin{Theorem}
\label{thm:exct:persist}
Let
\begin{align*}
c_1&=r_1-\dfrac{\sigma_1^2}{2},\\
c_2&=r_2-\dfrac{\sigma_2^2}{2},\\
c_4&=a_1K_1+a_2K_2-\delta -\dfrac{\sigma_3^2}{2}.
\end{align*}
For any given initial condition \eqref{ID}, if 
$$
\min\{c_1,c_2,1-r_1+2r_1/K_1,1-r_2+2r_2/K_2\}>0,
\quad \text{and} \quad 
c_4\leq 0, 
$$
then
\begin{align*}
\limsup_{t\rightarrow \infty}\left<z(t)\right>=0,
\end{align*}
and all the prey persist in mean. Moreover, we have
\begin{align*}
\liminf_{t\rightarrow \infty}\left<x(t)\right> 
& \dfrac{r_1-\frac{\sigma_1^2}{2}}{1-r_1+2r_1/K_1}\text{ a.s. }\\
\liminf_{t\rightarrow \infty}\left<y(t)\right> 
& \dfrac{r_2-\frac{\sigma_2^2}{2}}{1-r_2+2r_2/K_2}\text{ a.s. }
\end{align*}
\end{Theorem}

\begin{proof}
Using the first equation of \eqref{sys2}, we have
\begin{align*}
\dfrac{x(t)-x(0)}{t}
=&\dfrac{1}{t}\int_0^t\left(r_1x(s)\left(1-\dfrac{x(s-\tau_1)}{K_1}\right)
-\alpha_1 x(s)z(s)+\beta x(s)y(s)z(s)\right) ds\\ 
&+\sigma_1 x(s)dW_1(s)+ \dfrac{1}{t}
\int_0^t\left(\displaystyle\int_{U}q_1(u)x (s-)\tilde{N}(ds,du)\right)ds.
\end{align*}
Let us define
$$ 
F\left(x\right)=\log\left(x(t)\right)-\dfrac{r_1}{K_1}\int_t^{t+\tau_1}x(s-\tau_1)ds.
$$
Then, using It\^{o}'s formula, we have
\begin{align*}
dF=&r_1-\dfrac{\sigma_1^2}{2}+\int_U\log\left(1+q_1(u)\right)-q_1(u)\nu(du)\\
&-\dfrac{r_1x(t)}{K_1}-z(t)\left(\alpha_1-\beta y(t)\right)\; dt
+\sigma_1dW_1(t)+\int_U \log\left(1+q_1(u)\right)\tilde{N}(dt,du).
\end{align*}
Thus,
\begin{align*}
\dfrac{x(t)-x(0)}{t}+\dfrac{F(t)-F(0)}{t}
\geq & r_1\left<x\right>-\dfrac{r_1}{K_1}
\left<x\right>+r_1-\dfrac{\sigma_1^2}{2}+\alpha_1\dfrac{1}{t}
\int_0^t \int_U\log\left(1+q_1(u)\right)\\
&-q_1(u)\nu(du)ds-\dfrac{r_1x(t)}{K_1}
+\dfrac{1}{t}\int_0^t \sigma_1 (1+x(s))dW_1(s)\\
&+\dfrac{1}{t}\int_0^t\left(\displaystyle\int_{U}\log\left(1+q_1(u)\right)
+q_1(u)x (s-)\tilde{N}(ds,du)\right)ds.
\end{align*}
Using the strong law of large numbers for martingales, we obtain
$$
\liminf_{t\rightarrow +\infty} \left<x\right>
\geq \dfrac{r_1-\frac{\sigma_1^2}{2}}{1-r_1+2r_1/K_1}.
$$
With the same method, we prove that
$$
\liminf_{t\rightarrow +\infty} \left<y\right>
\geq \dfrac{r_2-\frac{\sigma_2^2}{2}}{1-r_2+2r_2/K_2}.
$$
Now, consider 
$$
V(t)=\log(z)+\int_t^{t+\tau_3}\left(a_1x(s-\tau_3)+a_2y(s-\tau_3)\right) ds.
$$
Using It\^{o}'s formula, we have that
\begin{equation*}
\begin{aligned}
dV=&LV\; dt+\sigma_3dW_3(t)
+\int_U \log\left(1+q_3(u)\right)\tilde{N}(dt,du),
\end{aligned}
\end{equation*}
where
\begin{align*}
LV=&-\delta-\alpha_3z+a_1x+a_2y-\dfrac{\sigma_3^2}{2}
+\int_U\log\left(1+q_3(u)\right)-q_3(u)\nu(du)\\
\leq &-\delta-\alpha_3z+a_1K_1+a_2K_2-\dfrac{\sigma_3^2}{2}.
\end{align*}
Then,
\begin{align*}
\dfrac{V(t)-V(0)}{t}
\leq & a_1K_1+a_2K_2-\delta -\alpha_3\left<z\right>-\dfrac{\sigma_3^2}{2}\\
&+\dfrac{1}{t}\int_0^t \left(\sigma_3dW_3(t)
+\int_U \log\left(1+q_3(u)\right)\tilde{N}(dt,du)\right),
\end{align*}
so that
\begin{align*}
\left<z\right>
\leq &\dfrac{1}{\alpha_3}\left(a_1K_1+a_2K_2-\delta 
-\dfrac{\sigma_3^2}{2}-\dfrac{V(t)-V(0)}{t}\right)\\
&+\dfrac{1}{\alpha_3}\dfrac{1}{t}\int_0^t \left(\sigma_3dW_3(t)
+\int_U \log\left(1+q_3(u)\right)\tilde{N}(dt,du)\right).
\end{align*}
We conclude that
$$
\limsup_{t\rightarrow +\infty}\left<z\right>  
\leq \dfrac{1}{\alpha_3}\left(a_1K_1+a_2K_2
-\delta -\dfrac{\sigma_3^2}{2}\right),
$$
and the proof is complete.
\end{proof}

% ----------------------------------

\section{Stochastic Persistence}
\label{sec:5}

This section is devoted to proving the persistence 
of the prey and predator populations.

\begin{Theorem}
\label{thm:persistence}
If
\begin{align*}
&\dfrac{r_1-\frac{\sigma_1^2}{2}}{1-r_1+2r_1/K_1}>0,\\
&\dfrac{r_2-\frac{\sigma_2^2}{2}}{1-r_2+2r_2/K_2}>0,\\
&a_1\left(\dfrac{r_1-\frac{\sigma_1^2}{2}}{1-r_1+2r_1/K_1}\right)
+a_2\left(\dfrac{r_2-\frac{\sigma_2^2}{2}}{1-r_2+2r_2/K_2}\right)
-\delta -\dfrac{\sigma_3^2}{2}>0,
\end{align*} 
and
$$
\min\{1-r_1+2r_1/K_1,1-r_2+2r_2/K_2\}> 0,
$$
then all the prey and predators persist in mean. Moreover, one has
\begin{gather*}
\liminf_{t\rightarrow +\infty} \left<x\right>
\geq \dfrac{r_1-\frac{\sigma_1^2}{2}}{1-r_1+2r_1/K_1}\text{ a.s. }\\
\liminf_{t\rightarrow +\infty} \left<y\right>
\geq \dfrac{r_2-\frac{\sigma_2^2}{2}}{1-r_2+2r_2/K_2} \text{ a.s. }\\
\liminf_{\rightarrow +\infty}\left<z\right>  
\geq \dfrac{1}{\alpha_3}\left(a_1\left(\dfrac{r_1
-\frac{\sigma_1^2}{2}}{1-r_1+2r_1/K_1}\right)
+a_2\left(\dfrac{r_2-\frac{\sigma_2^2}{2}}{1-r_2+2r_2/K_2}\right)
-\delta -\dfrac{\sigma_3^2}{2}\right) \text{ a.s. }
\end{gather*}
\end{Theorem}

\begin{proof}
The first equation of \eqref{sys2} implies that
\begin{align*}
\dfrac{x(t)-x(0)}{t}
=&\dfrac{1}{t}\int_0^t\left(r_1x(s)\left(1-\dfrac{x(s-\tau_1)}{K_1}\right)
-\alpha_1 x(s)z(s)+\beta x(s)y(s)z(s)\right) ds \\
&+\sigma_1 x(s)dW_1(s)+ \dfrac{1}{t}
\int_0^t\left(\displaystyle\int_{U}q_1(u)x (s-)\tilde{N}(ds,du)\right)ds.
\end{align*}
Let 
$$ 
F\left(x\right)=\log\left(x(t)\right)-\dfrac{r_1}{K_1}\int_t^{t+\tau_1}x(s-\tau_1)ds.
$$
Then, from It\^{o}'s formula,
\begin{align*}
dF=
&r_1-\dfrac{\sigma_1^2}{2}
+\int_U\log\left(1+q_1(u)\right)-q_1(u)\nu(du)\\
&-\dfrac{r_1x(t)}{K_1}-z(t)\left(\alpha_1-\beta y(t)\right)\; dt
+\sigma_1dW_1(t)+\int_U \log\left(1+q_1(u)\right)\tilde{N}(dt,du).
\end{align*}
Therefore,
\begin{align*}
\dfrac{x(t)-x(0)}{t}+\dfrac{F(t)-F(0)}{t}
\geq & r_1\left<x\right>-\dfrac{r_1}{K_1}\left<x\right>
-\alpha_1\dfrac{1}{t}\int_0^t(x+1)zds\\
&+r_1-\dfrac{\sigma_1^2}{2}+\alpha_1\dfrac{1}{t}\int_0^t 
\int_U\log\left(1+q_1(u)\right)-q_1(u)\nu(du)ds\\
&-\dfrac{r_1x(t)}{K_1}+\dfrac{1}{t}\int_0^t \sigma_1 (1+x(s))dW_1(s)\\
&+\dfrac{1}{t}\int_0^t\left(\displaystyle\int_{U}\log\left(1+q_1(u)\right)
+q_1(u)x (s-)\tilde{N}(ds,du)\right)ds,
\end{align*}
and, using the strong law of large numbers for martingales, we obtain
$$
\liminf_{\rightarrow +\infty} \left<x\right>
\geq \dfrac{r_1-\frac{\sigma_1^2}{2}}{1-r_1+2r_1/K_1}.
$$
Similarly, we prove that
$$
\liminf_{t\rightarrow +\infty} \left<y\right>
\geq \dfrac{r_2-\frac{\sigma_2^2}{2}}{1-r_2+2r_2/K_2}.
$$
Using now It\^{o}'s formula with
$$
V(t)=\log(z)+\int_t^{t+\tau_3}\left(a_1x(s-\tau_3)+a_2y(s-\tau_3)\right) ds,
$$
we have 
\begin{equation*}
\begin{aligned}
dV=&LV\; dt+\sigma_3dW_3(t)
+\int_U \log\left(1+q_3(u)\right)\tilde{N}(dt,du),
\end{aligned}
\end{equation*}
where
\begin{align*}
LV=&-\delta-\alpha_3z+a_1x+a_2y-\dfrac{\sigma_3^2}{2}
+\int_U\log\left(1+q_3(u)\right)-q_3(u)\nu(du).
\end{align*}
Then,
\begin{align*}
\dfrac{V(t)-V(0)}{t}
\geq & a_1\left(\dfrac{r_1-\frac{\sigma_1^2}{2}}{1-r_1+2r_1/K_1}\right)
+a_2\left(\dfrac{r_2-\frac{\sigma_2^2}{2}}{1-r_2+2r_2/K_2}\right)-\delta 
-\alpha_3\left<z\right>-\dfrac{\sigma_3^2}{2}\\
&+\dfrac{1}{t}\int_0^t \left(\sigma_3dW_3(t)
+\int_U \log\left(1+q_3(u)\right)\tilde{N}(dt,du)\right)
\end{align*}
so that
\begin{align*}
\left<z\right>
\geq &\dfrac{1}{\alpha_3}\left(a_1\left(
\dfrac{r_1-\frac{\sigma_1^2}{2}}{1-r_1+2r_1/K_1}\right)
+a_2\left(\dfrac{r_2-\frac{\sigma_2^2}{2}}{1-r_2+2r_2/K_2}\right)
-\delta -\dfrac{\sigma_3^2}{2}-\dfrac{V(t)-V(0)}{t}\right)\\
&+\dfrac{1}{\alpha_3}\dfrac{1}{t}\int_0^t \left(\sigma_3dW_3(t)
+\int_U \log\left(1+q_3(u)\right)\tilde{N}(dt,du)\right).
\end{align*}
Therefore,
$$
\liminf_{\rightarrow +\infty}\left<z\right>  
\geq \dfrac{1}{\alpha_3}\left(a_1\left(\dfrac{r_1
-\frac{\sigma_1^2}{2}}{1-r_1+2r_1/K_1}\right)
+a_2\left(\dfrac{r_2-\frac{\sigma_2^2}{2}}{1-r_2+2r_2/K_2}\right)
-\delta -\dfrac{\sigma_3^2}{2}\right),
$$
which proves the intended result.
\end{proof}

% ----------------------------------

\section{Numerical Simulations}
\label{sec:6}

\textls[-15]{This section is devoted to illustrating our mathematical findings using numerical simulations. }

In the following examples, we apply the algorithm presented in \cite{33} 
to solve system~\eqref{sys2}. In our simulations, the time 
period is in days. The different values for the parameters 
used in our numerical simulations are given in Table~\ref{tabl1}. 

% -----------------------------------
\begin{table}[H]
\caption{Values of the parameters used in the numerical simulations.}
\label{tabl1}
\newcolumntype{C}{>{\centering\arraybackslash}X}
\begin{tabularx}{\textwidth}{CCCC}
\toprule 
\textbf{Parameters} & \textbf{Figure~\ref{fig1}} 
& \textbf{Figure~\ref{fig2}} & \textbf{Figure~\ref{fig3}}\\ \midrule
$r_1$ 	& $0.7$& $1.7$ & $2$\\
$K_1$	& $100$   &  $100$& $100$\\
$r_2$ & $0.65$&$1.8$&$2.3$ \\
$K_2$ & $100 $   & $100$& $100$ \\
$\alpha_1$& $0.3$  & $0.2$&$0.13$\\
$\alpha_2$& $0.35$ &$0.28$&$0.17$\\
$\alpha_3$& $0.5$ &$0.5$&$0.2$\\
$\beta$& $0.0001$ &$0.0001$&$0.001$\\
$\delta$& $0.1$ &$0.4$&$0.02$\\
$\tau_1$& $0.5$ &$0.5$&$0.5$\\
$\tau_2$& $1$ &$1$&$1$\\
$\tau_3$& $1.5$ &$1.5$&$1.5$\\
$\sigma_1$ & $10^{-4}$&$10^{-5}$&$10^{-5}$ \\
$\sigma_2$ & $2\times 10^{-4}$&$2\times 10^{-4}$&$2\times 10^{-4}$\\
$\sigma_3$ & $2\times 10^{-4}$&$2\times 10^{-3}$&$2\times 10^{-3}$ \\
$q_1(u)$& $-0.04$&$-0.04$& $-0.04$ \\
$q_2(u)$& $-0.006$&$-0.006$& $-0.006$\\
$q_3(u)$& $-0.008$&$-0.008$&$-0.008$\\ \bottomrule
\end{tabularx}
\end{table}
% -----------------------------------

\begin{figure}[H]
\hspace{-30pt}\includegraphics[scale=0.4]{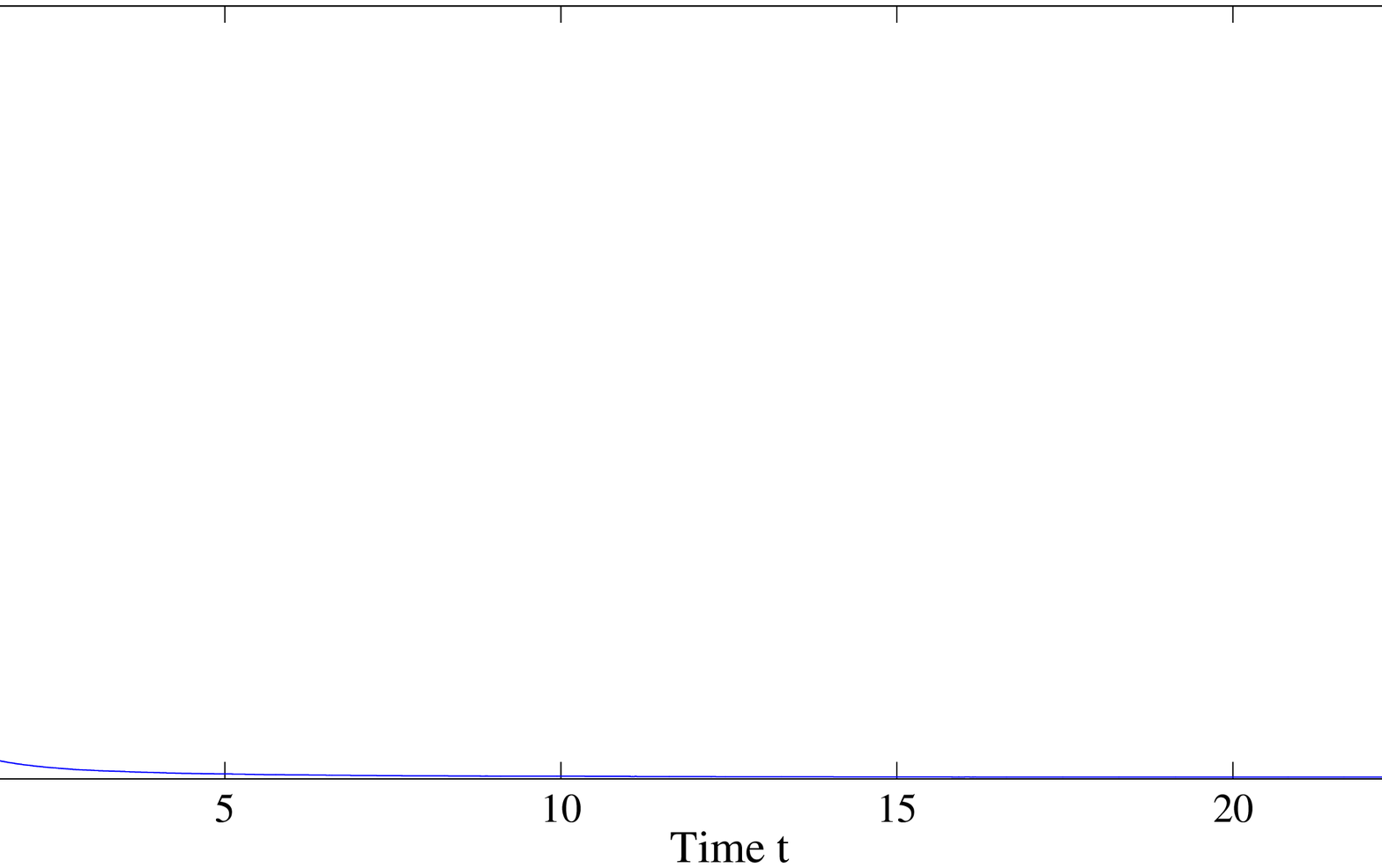}\\

\hspace{-30pt}\includegraphics[scale=0.4]{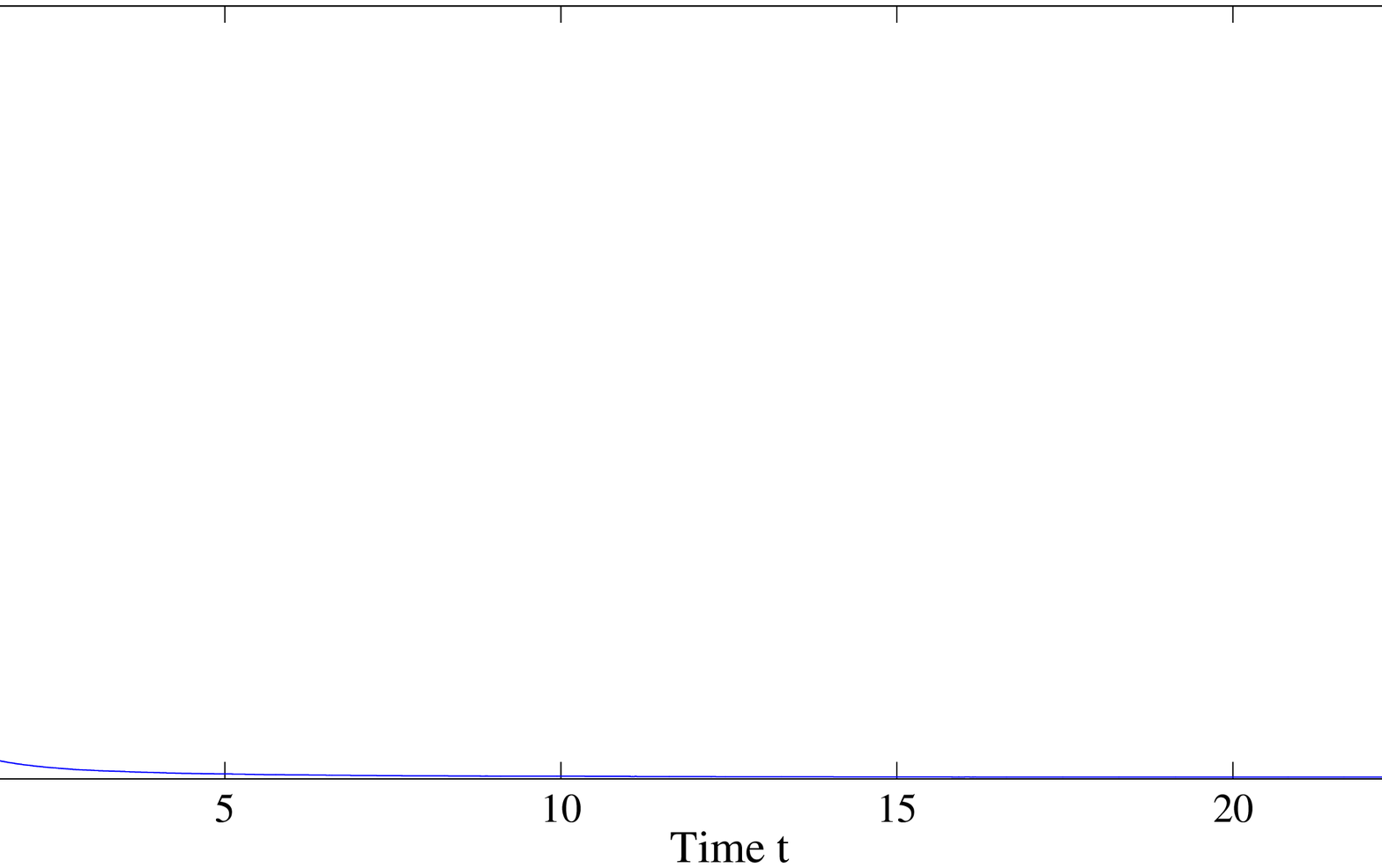}\\

\hspace{-30pt}\includegraphics[scale=0.4]{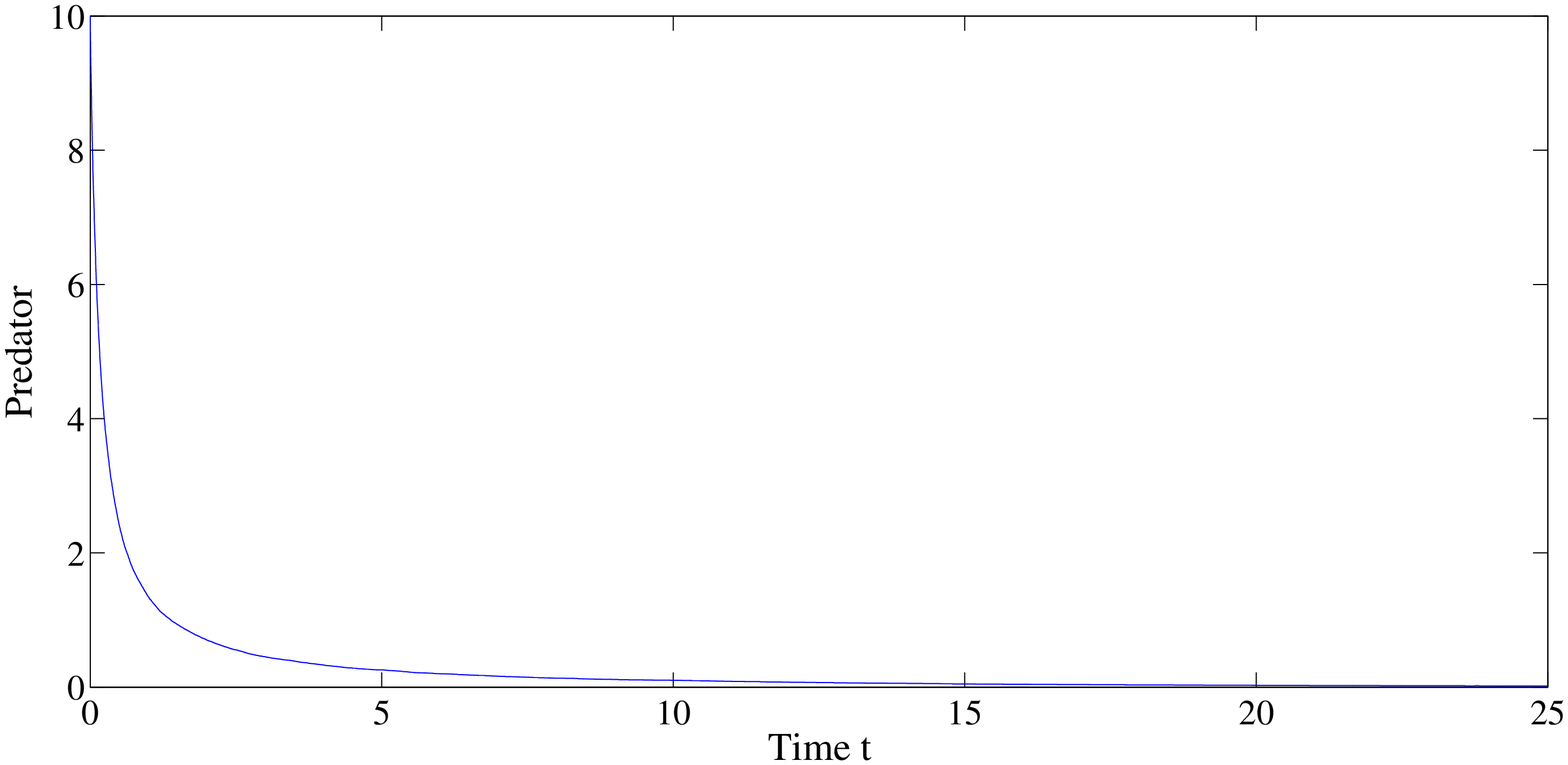}
\caption{The behavior of the prey and predator populations 
in the extinction situation described by 
Theorem~\ref{thm:extinction:prey:predator}.}
\label{fig1}
\end{figure}

Figure~\ref{fig1} shows the evolution of the two prey, $x$ and $y$, and the predator $z$,
in the case that both prey and predator vanish. This means the extinction of the populations,
which is in agreement with Section~\ref{sec:3}.

\begin{figure}[H]

\hspace{-27pt}\includegraphics[scale=0.4]{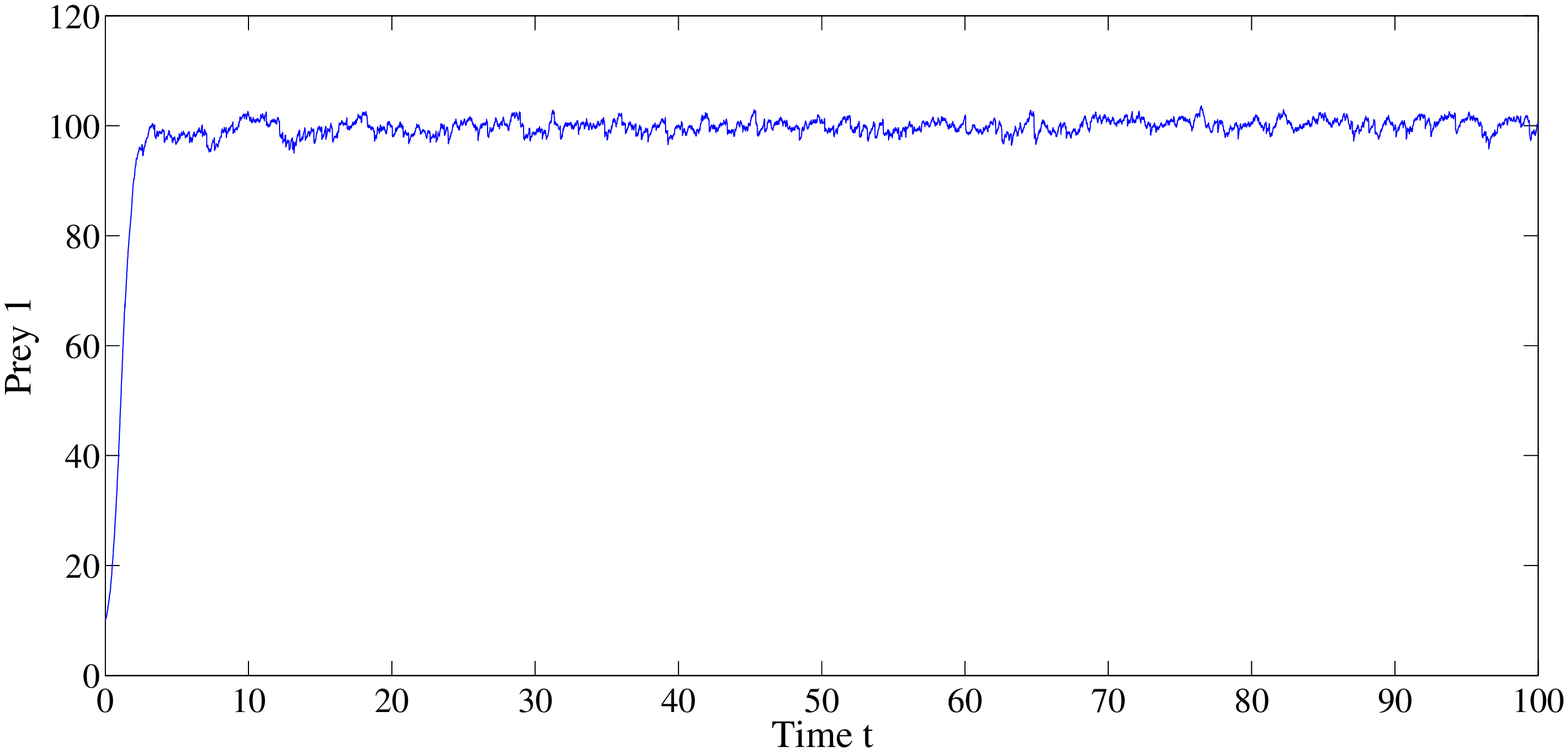}\\

\hspace{-27pt}\includegraphics[scale=0.4]{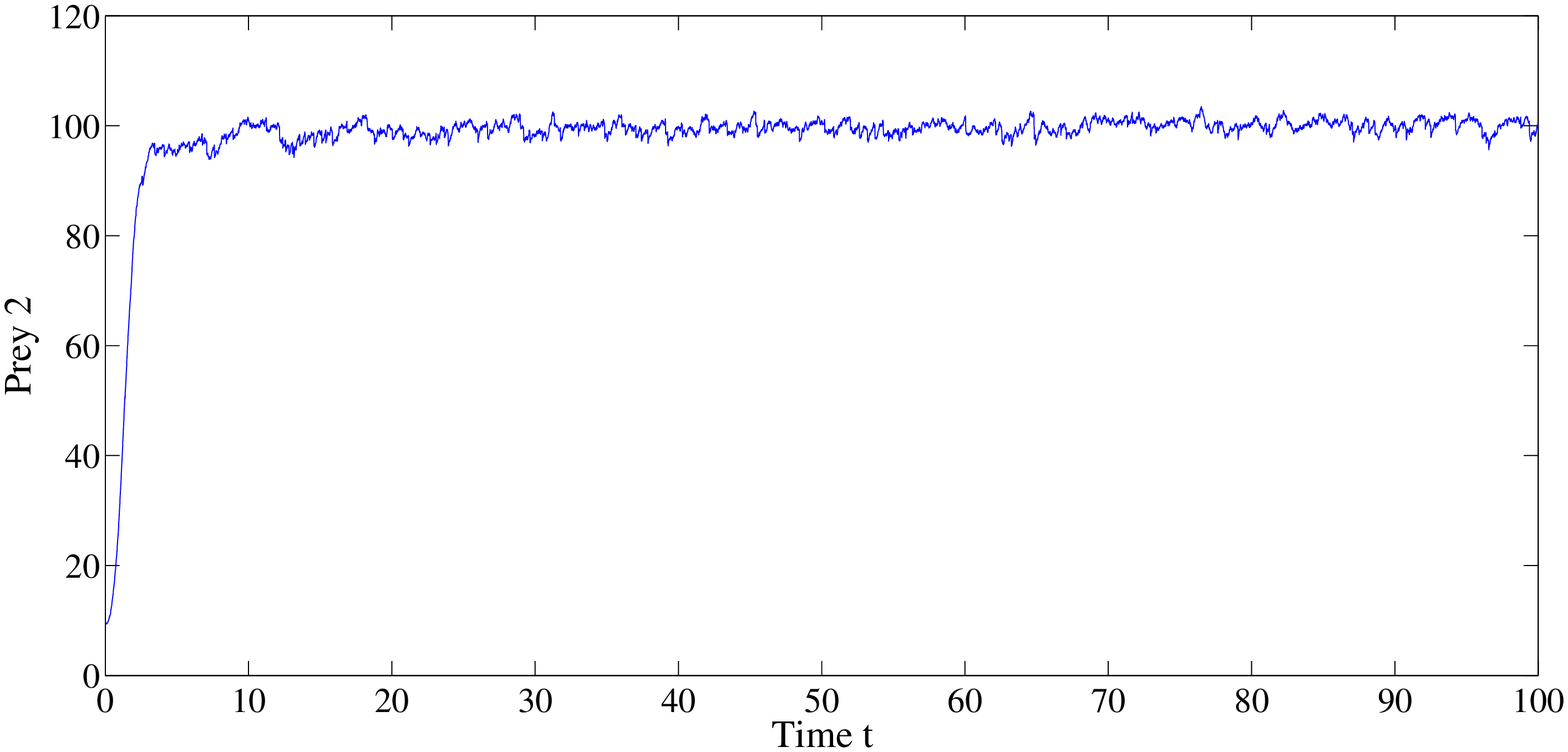}\\

\hspace{-27pt}\includegraphics[scale=0.4]{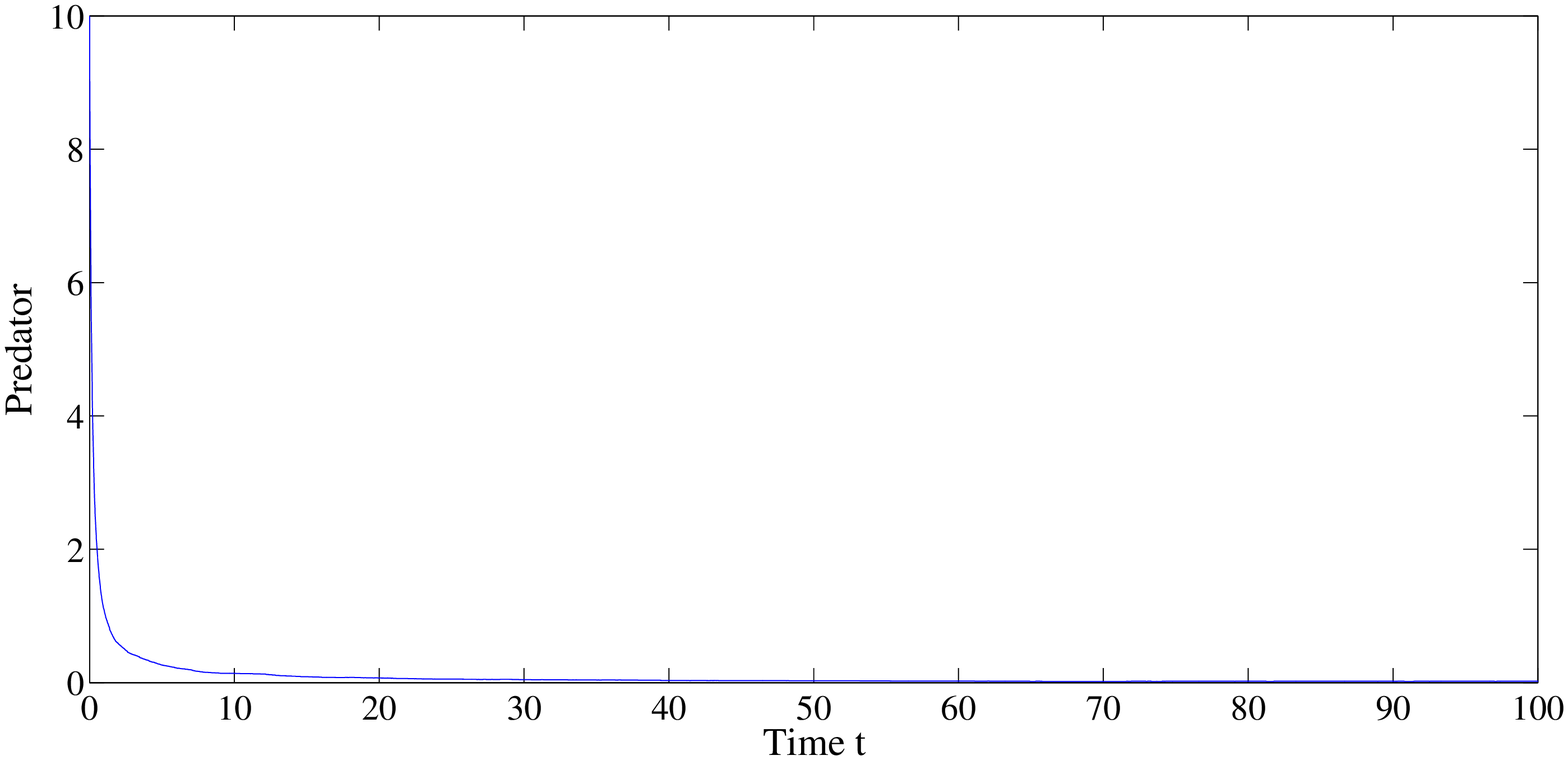}
\caption{The behavior of the prey and predator populations 
in the extinction--persistence situation described 
by Theorem~\ref{thm:exct:persist}.}
\label{fig2}
\end{figure}

Figure~\ref{fig2} presents the dynamics of $x$, $y$ and $z$ in the case 
studied in Section~\ref{sec:4}, where the predator tends to zero 
and the prey remains strictly positive: the population of predators will be extinct 
while the population of prey persists.

\begin{figure}[H]

\hspace{-30pt}\includegraphics[scale=0.4]{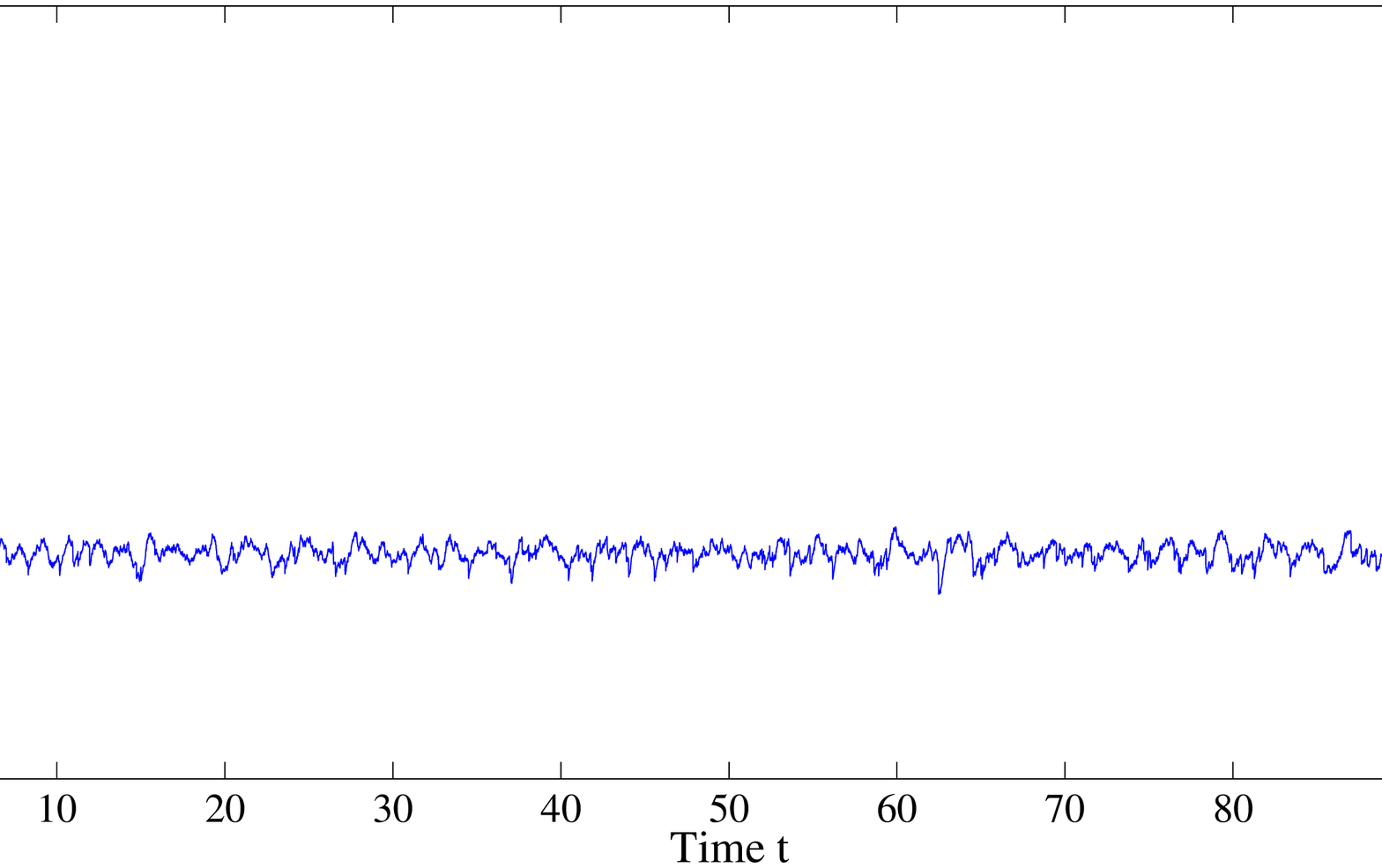}\\

\hspace{-30pt}\includegraphics[scale=0.4]{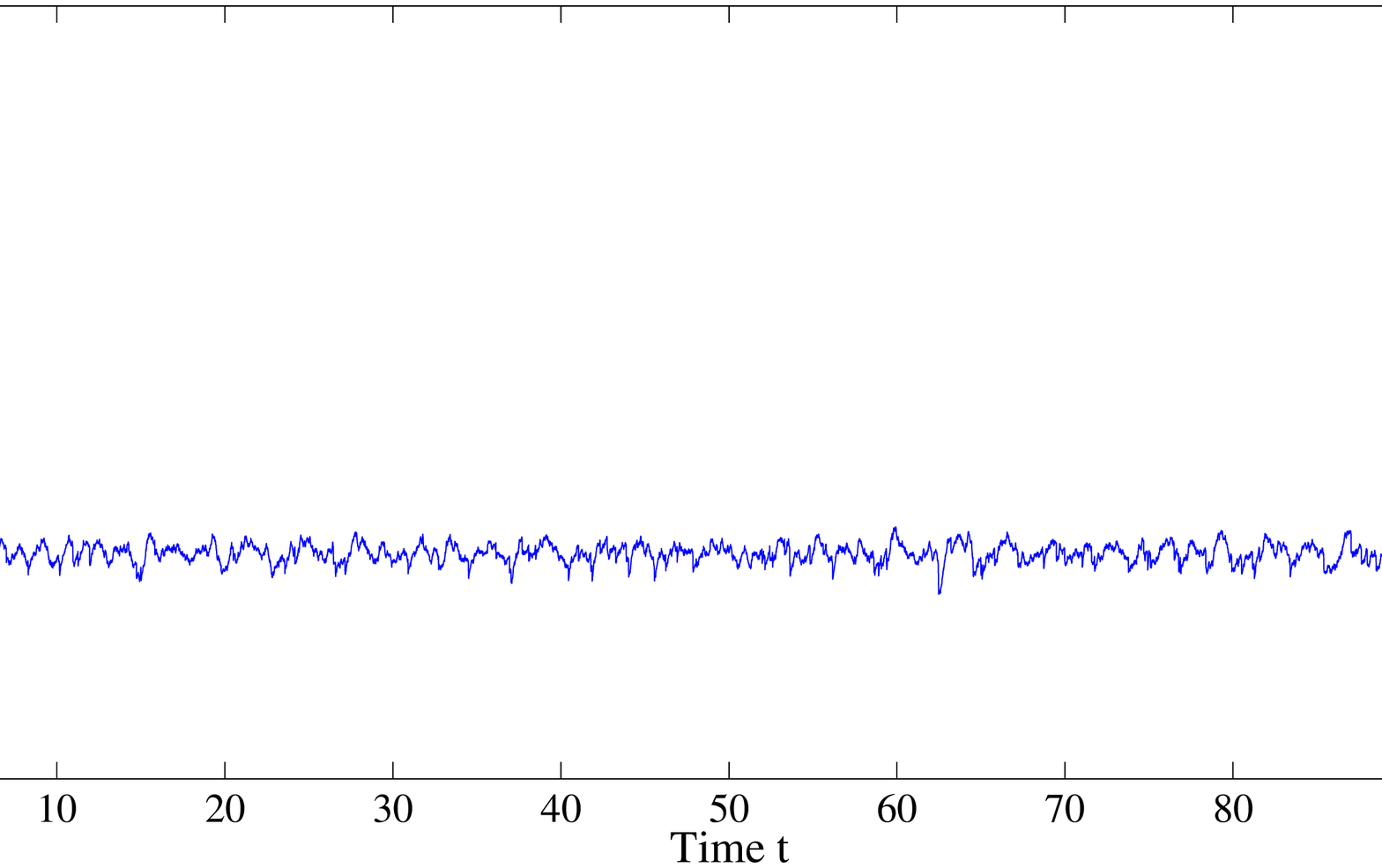}\\

\hspace{-30pt}\includegraphics[scale=0.4]{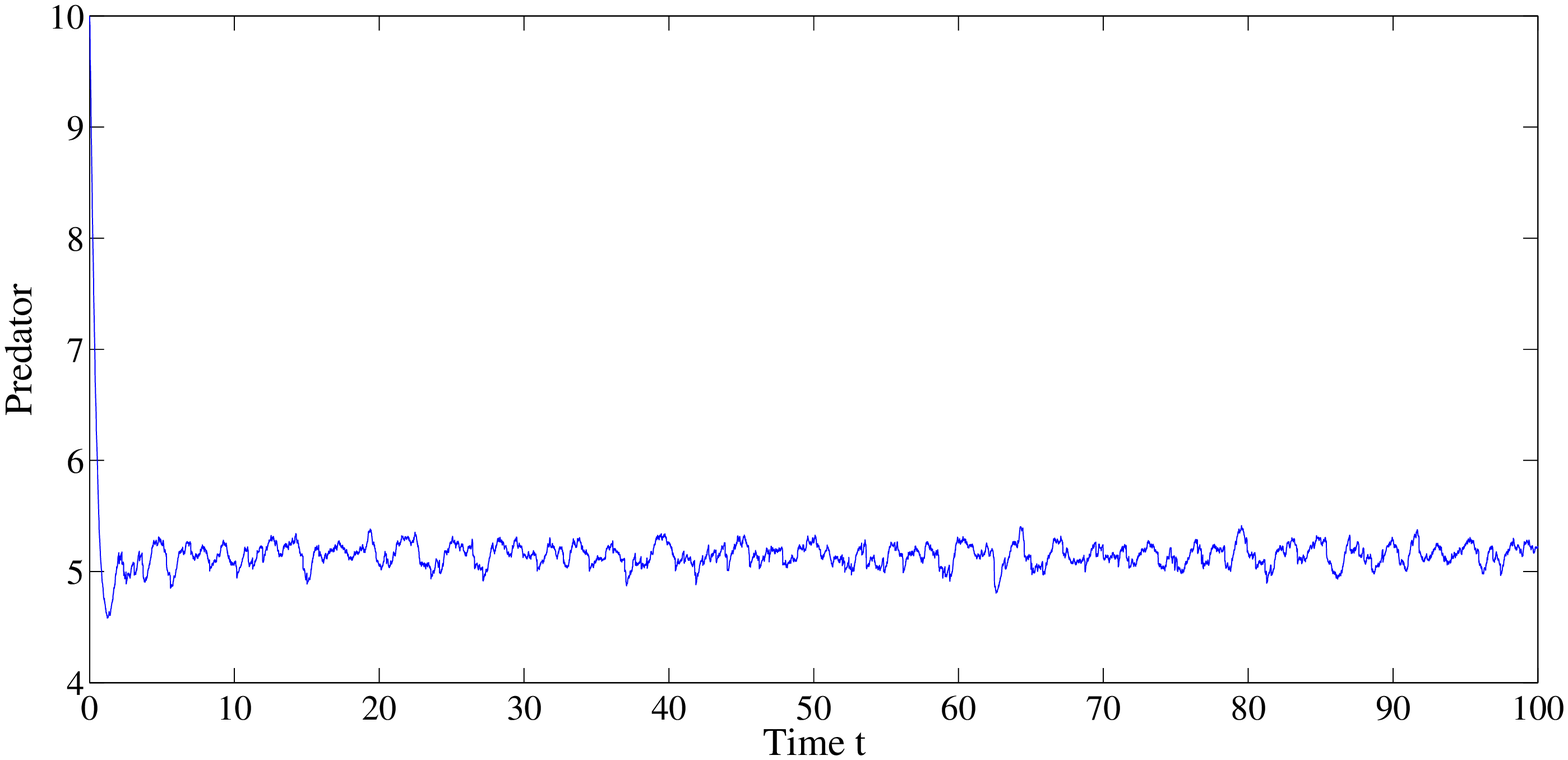}
\caption{The behavior of the prey and predator populations in the 
case of persistence described by Theorem~\ref{thm:persistence}.}
\label{fig3}
\end{figure}

Figure~\ref{fig3} illustrates the behavior of $x$, $y$ and $z$
in a situation where both the predator and the preys remain strictly positive,
i.e., all populations persist. This agrees 
with our theoretical result of Section~\ref{sec:5}.

\begin{figure}[H]
\hspace{-30pt}\includegraphics[scale=0.37]{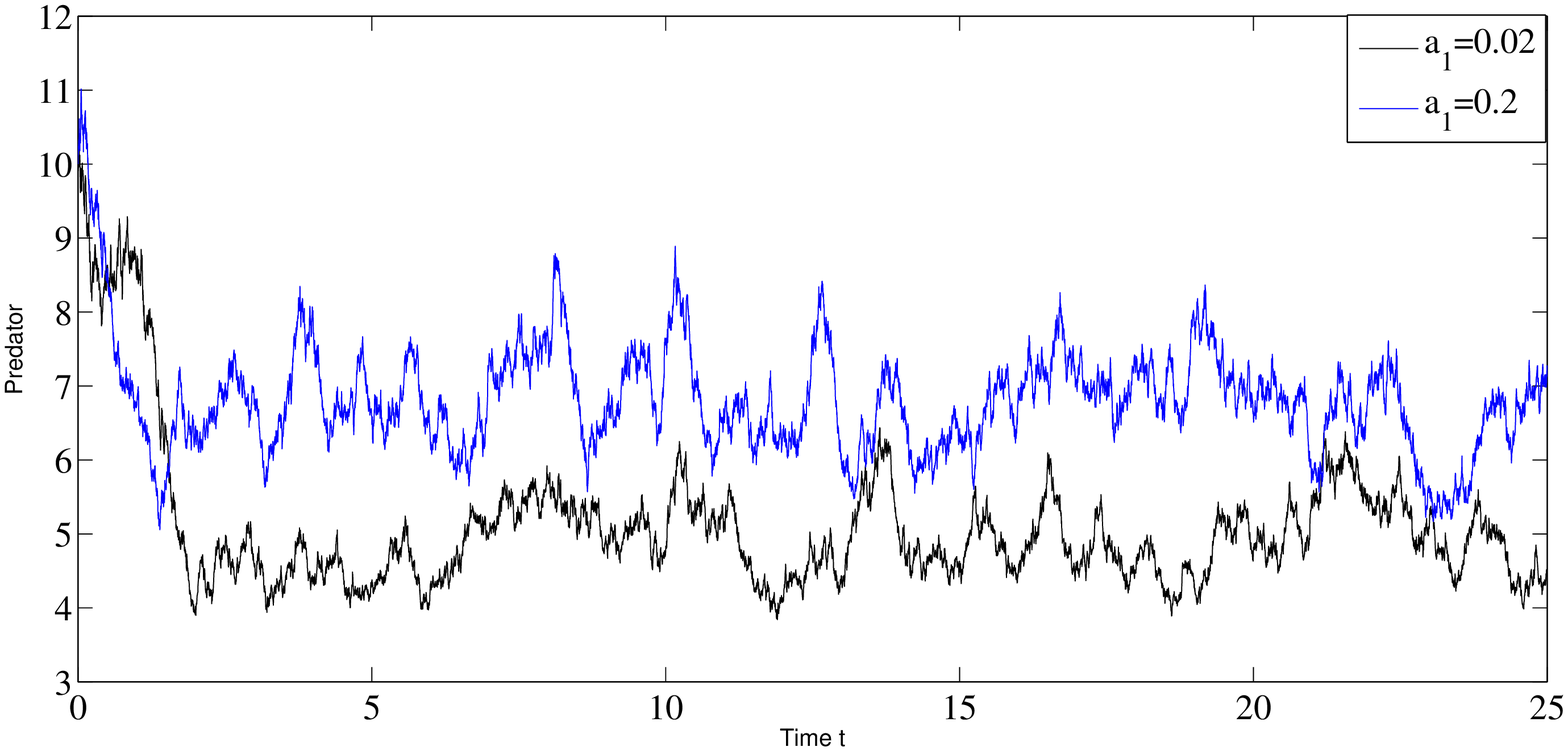}\\

\hspace{-30pt}\includegraphics[scale=0.37]{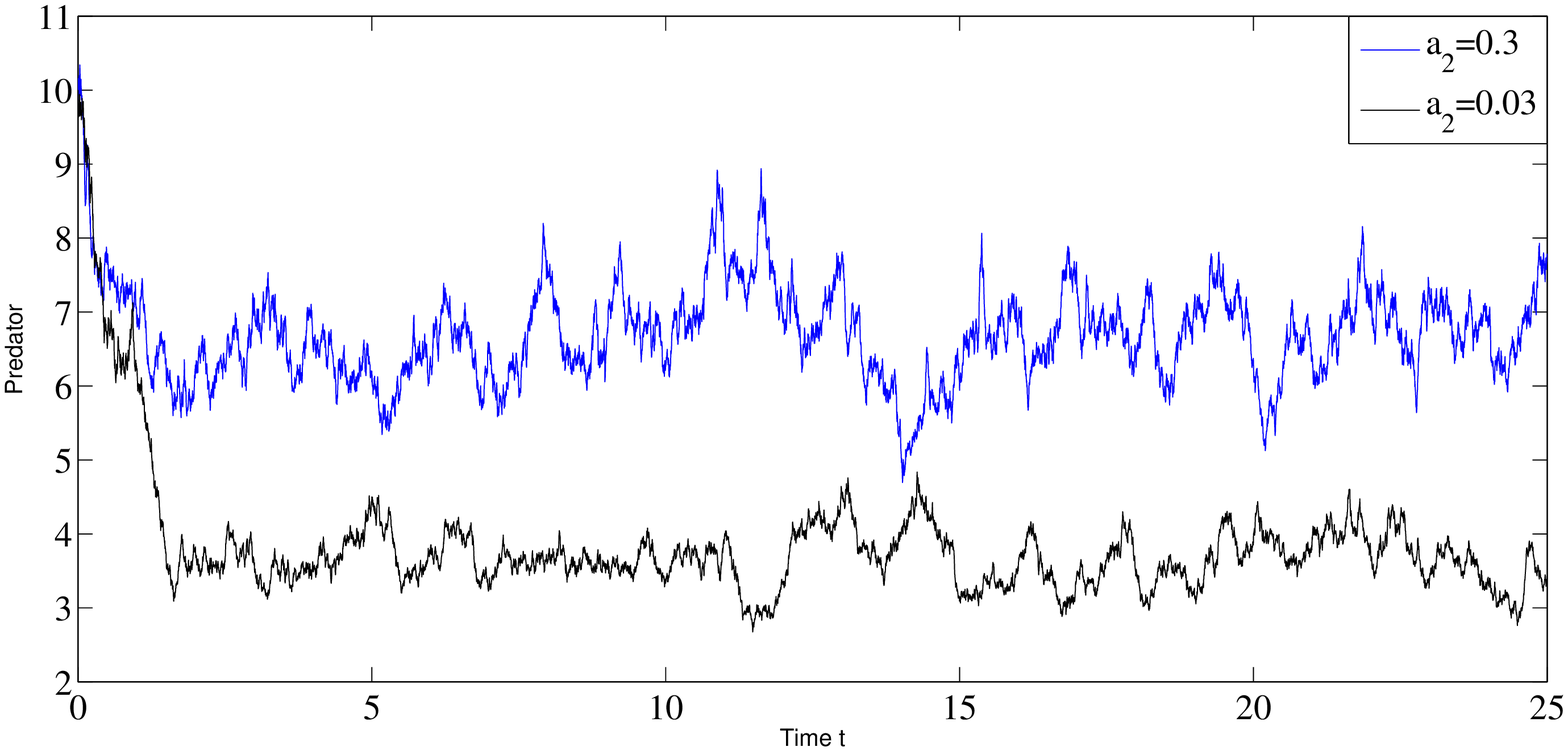}
\caption{The behavior of the predator population 
for different values of $a_1$ and $a_2$.}
\label{fig4}
\end{figure}

Figure~\ref{fig4} shows the impact of 
the transformation rates of predator to prey, $a_1$ and $a_2$, 
based on the dynamics of the predator. We can remark that $a_2$ 
has more influence on the number of predators than the parameter $a_1$.

\begin{figure}[H]
\hspace{-30pt}\includegraphics[scale=0.37]{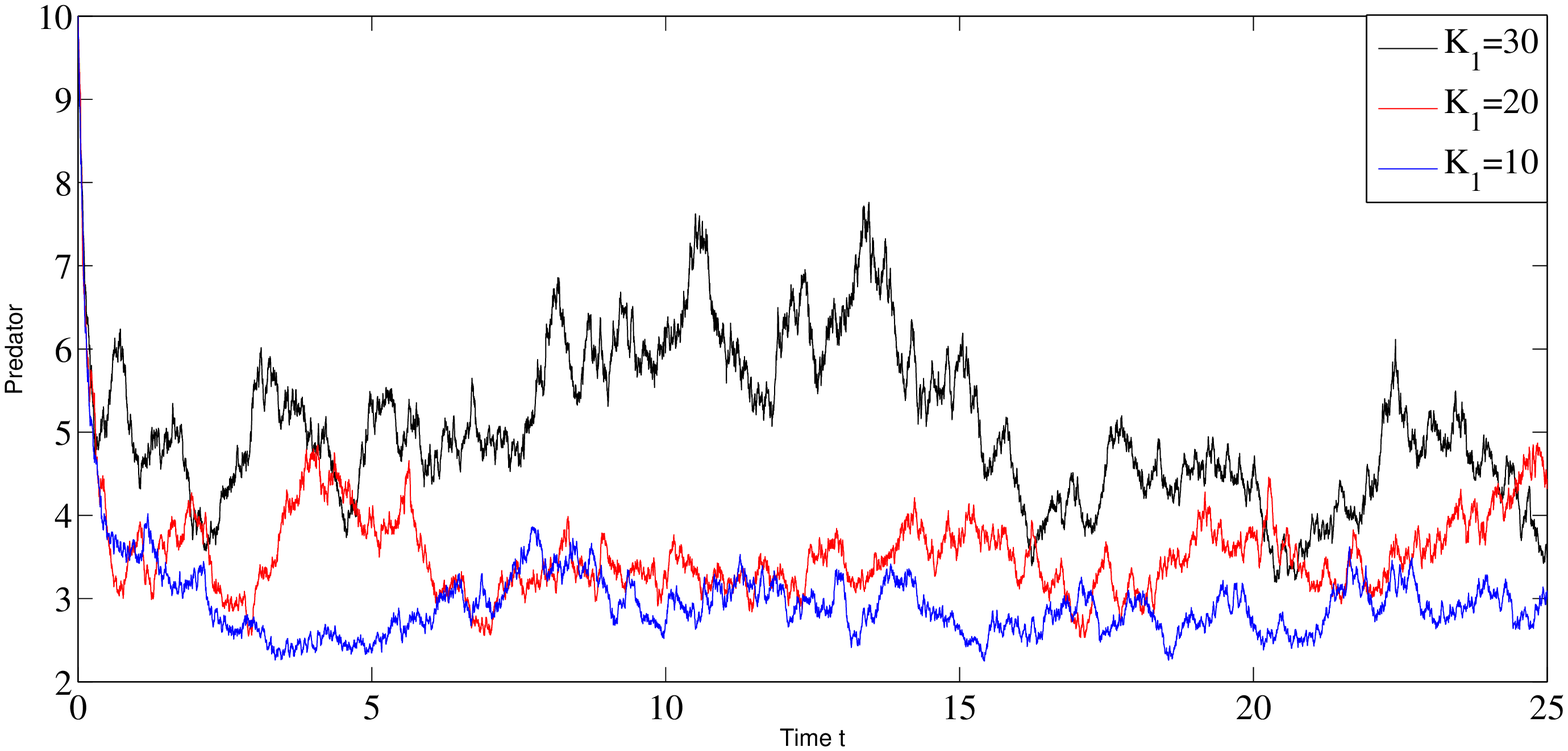}\\

\hspace{-30pt}\includegraphics[scale=0.37]{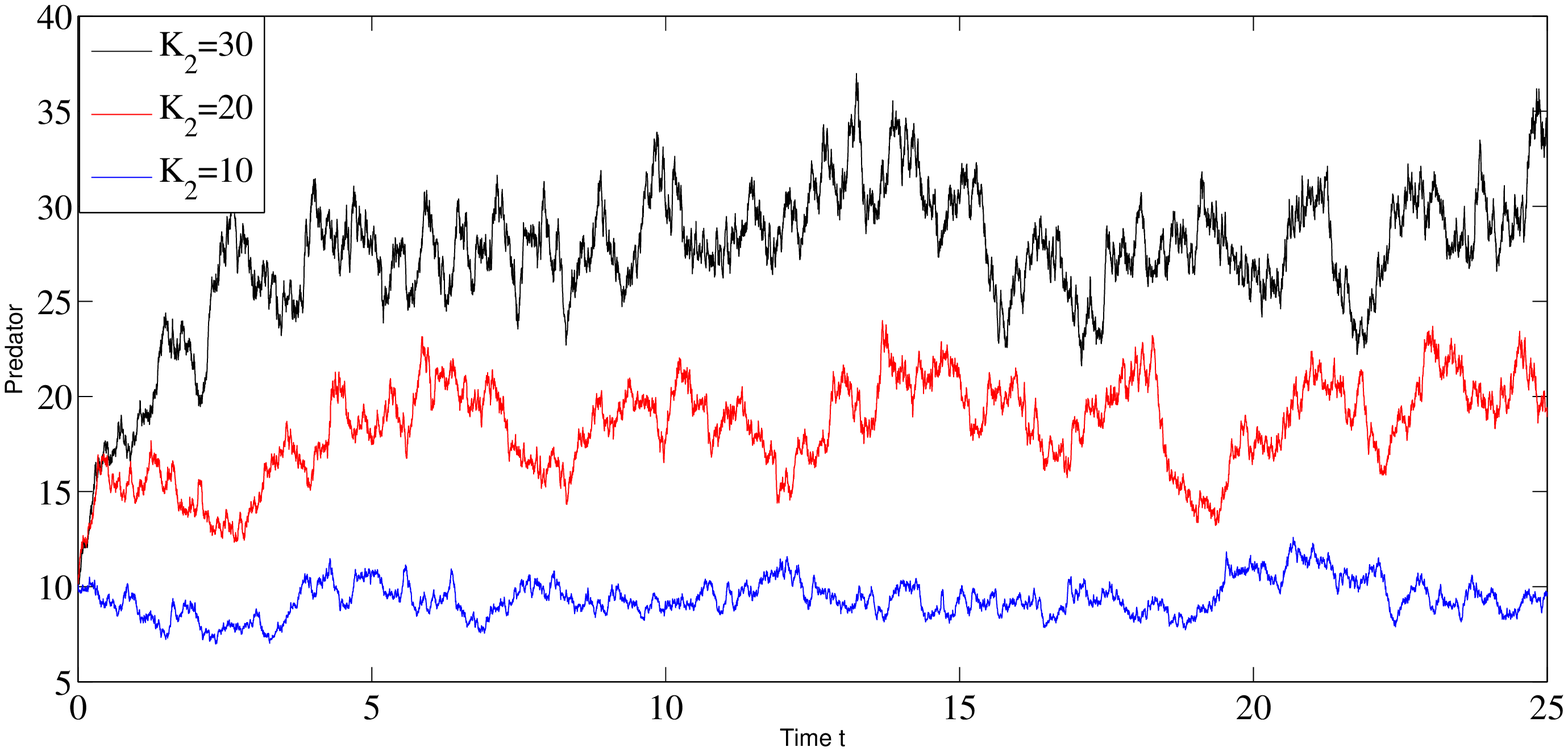}
\caption{The behavior~of the predator population 
for different values of $K_1$ and $K_2$.}
\label{fig5}
\end{figure}

Figure~\ref{fig5} shows the impact of the carrying capacities 
$K_1$ and $K_2$ on the dynamics of the predator. We observe that 
$K_2$ has more influence on the number of predators than the parameter $K_1$.

{It is known that delays can affect the behavior of the studied 
population \cite{Round2:01,Round2:02}. For this reason, now we study 
the effect of the delay on prey $1$, prey $2$, and the predator.}

% -----------------------------------
\unskip

% -----------------------------------
% -----------------------------------
\unskip

% -----------------------------------
% -----------------------------------
\unskip

\unskip
% -----------------------------------

\begin{figure}[H]
\hspace{-32pt}\includegraphics[scale=0.4]{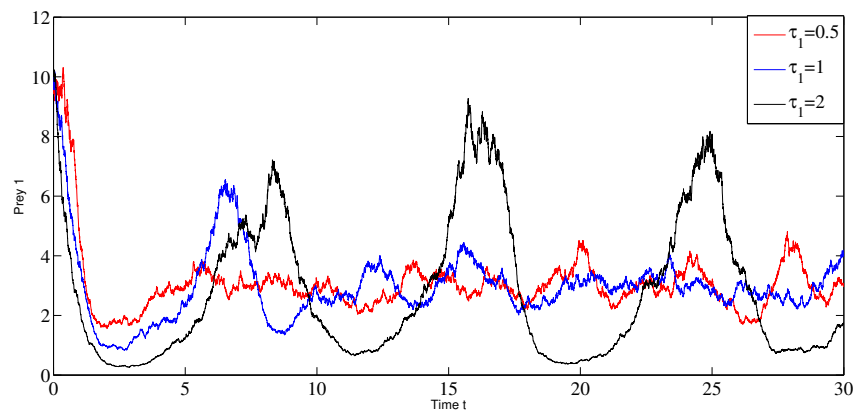}
\caption{The behavior of the prey 1 population
for different values of $\tau_1$.}
\label{fig6}
\end{figure}

{Figure~\ref{fig6} represents the dynamics of prey $1$ 
for different values of $\tau_1$. We observe that when one increases 
the value of $\tau_1$ from $0.5$ to $2$ days, that increases the amplitude 
of the prey $1$ with a delay translation.}

\vspace{-8pt}

% -----------------------------------

\begin{figure}[H]
\hspace{-32pt}\includegraphics[scale=0.4]{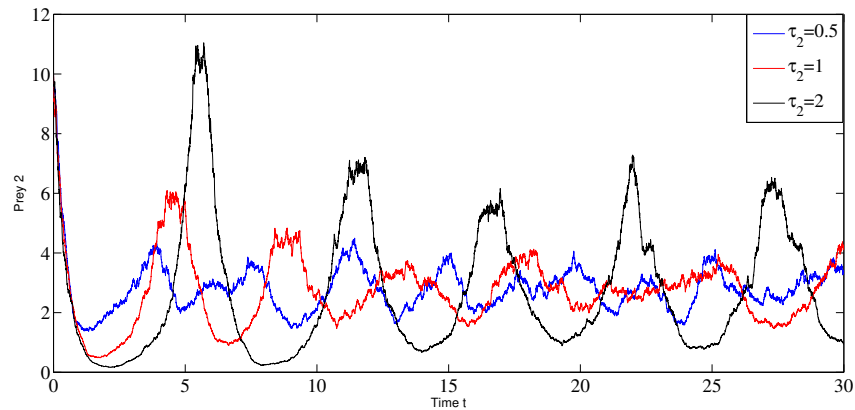}
\caption{The behavior of the prey 2 population 
for different values of $\tau_2$.}
\label{fig7}
\end{figure}

{Figure~\ref{fig7} represents the dynamics of prey $2$ for different 
values of $\tau_2$. We conclude that if we increase the value of $\tau_2$ 
from $0.5$ to $2$ days, then the amplitude of the prey $2$ 
increases with a delayed translation.}
% -----------------------------------
\vspace{-8pt}
% -----------------------------------

\begin{figure}[H]
\hspace{-32pt}\includegraphics[scale=0.4]{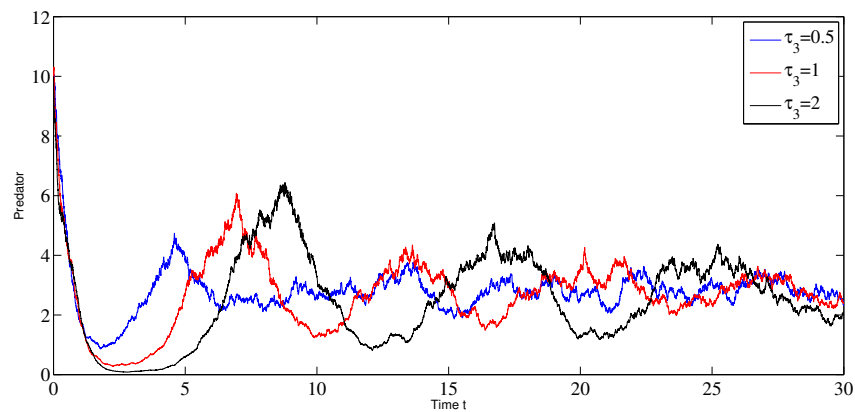}
\caption{The behavior of the predator population 
for different values of $\tau_3$.}
\label{fig8}
\end{figure}

{In Figure~\ref{fig8}, we show the behavior of the predator population for different 
values of $\tau_3$. We note that an increase of the value of $\tau_3$ from $0.5$ to $2$ days 
results in an increase of the amplitude of the predators with a delayed response.}
% -----------------------------------
\vspace{-8pt}
% -----------------------------------

\begin{figure}[H]
\hspace{-30pt}\includegraphics[scale=0.35]{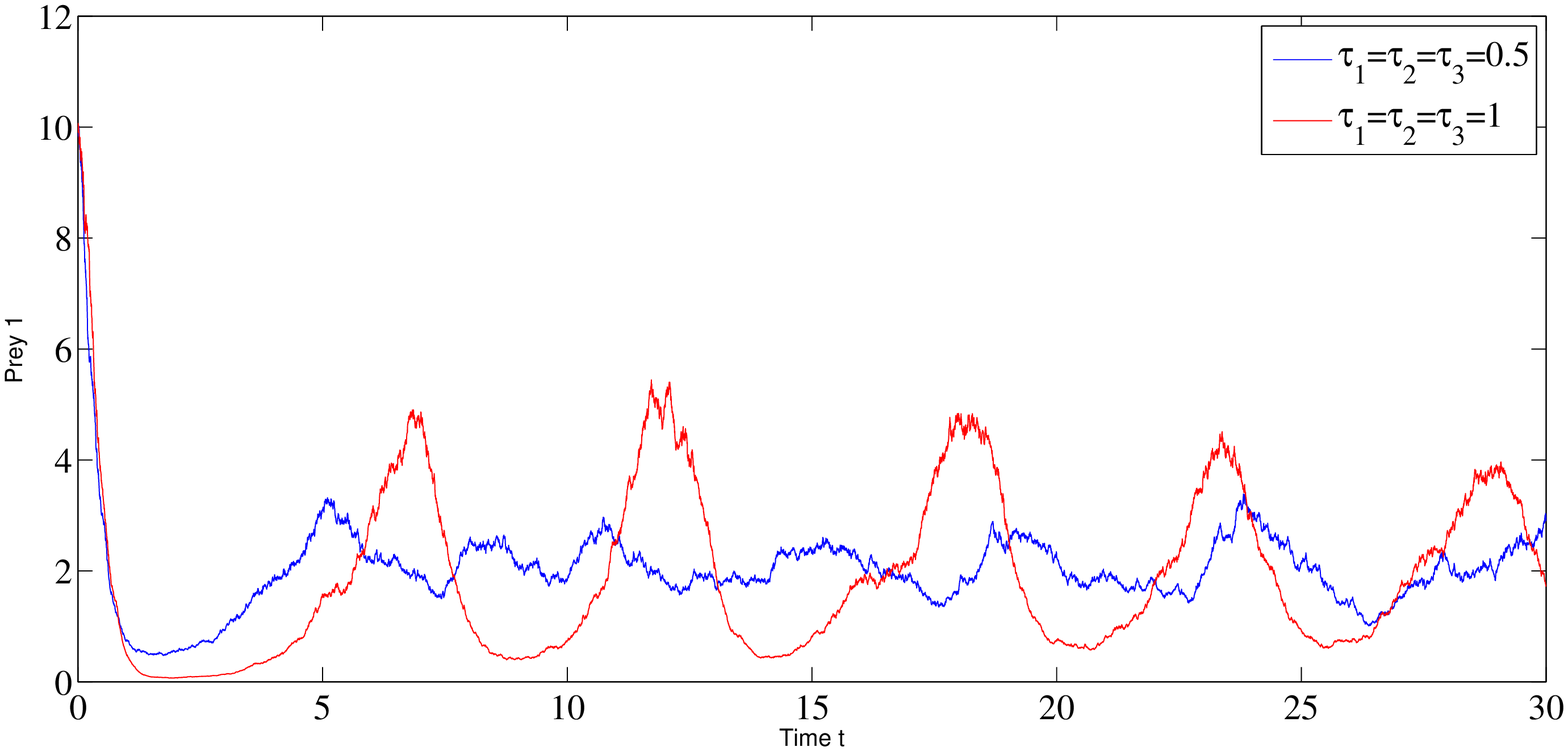}\\

\hspace{-30pt}\includegraphics[scale=0.35]{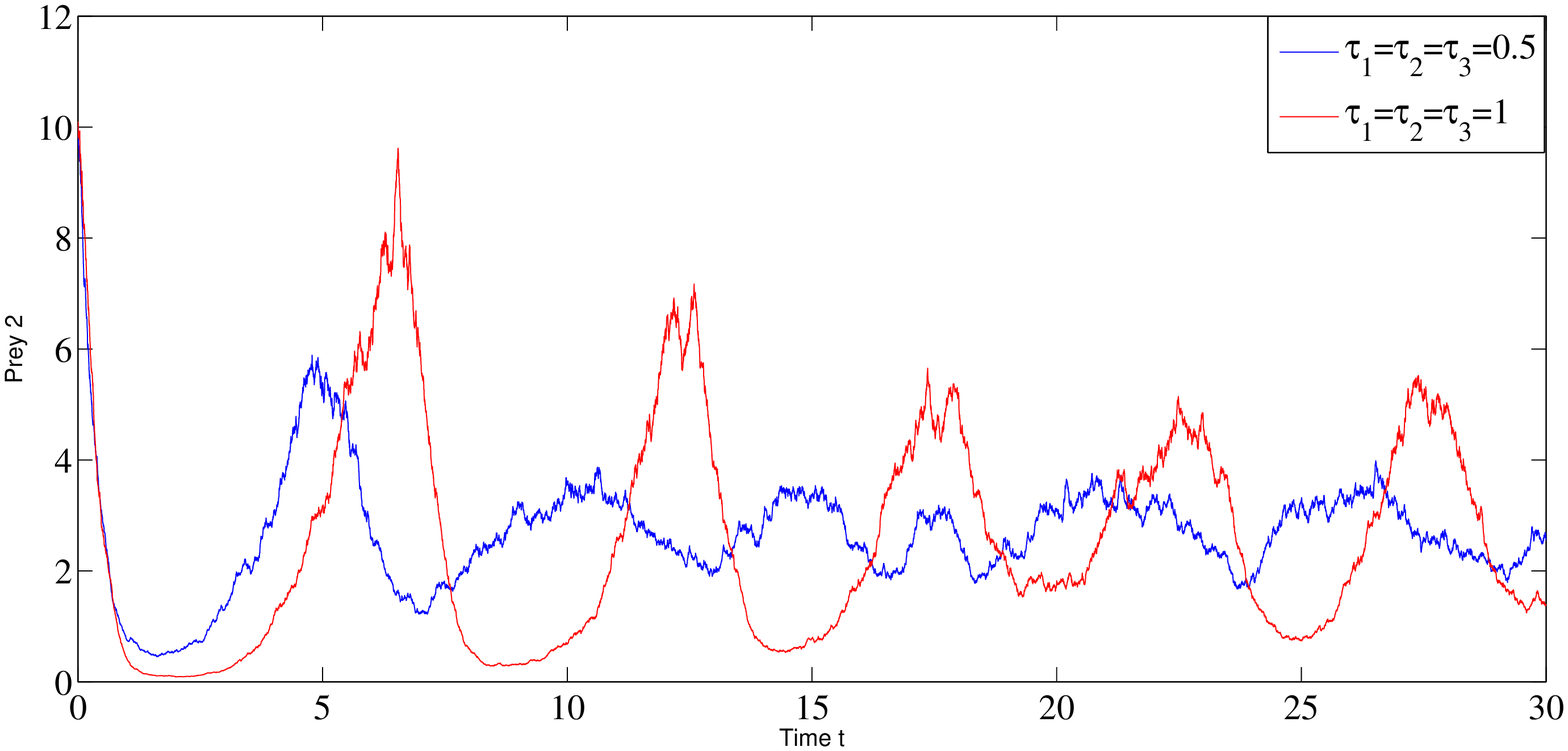}\\

\hspace{-30pt}\includegraphics[scale=0.35]{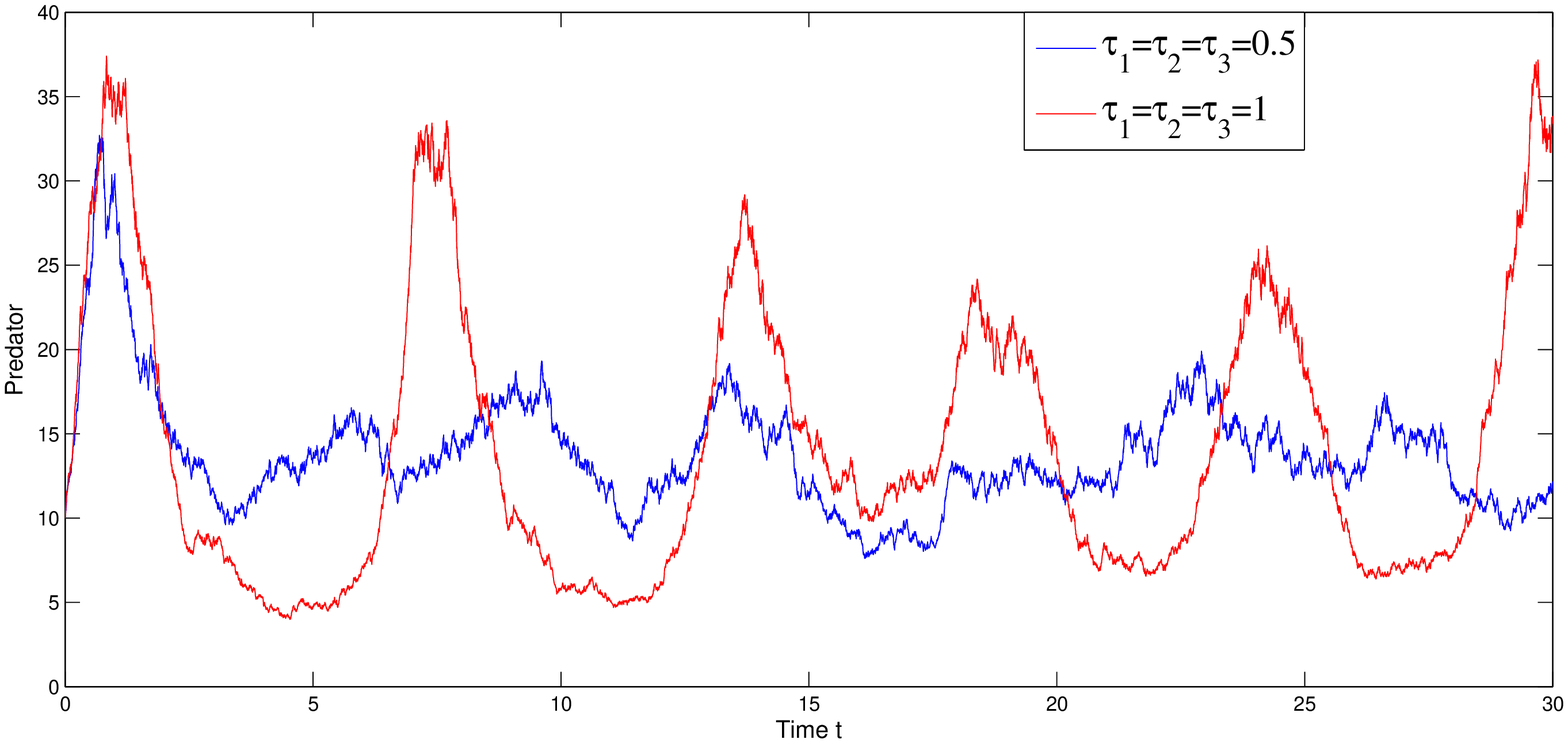}
\caption{The behavior of the populations 
for different values of $\tau_1$, $\tau_2$ and $\tau_3$.}
\label{fig9}
\end{figure}

{Finally, in Figure~\ref{fig9} we illustrate the effect of the three delays. 
We observe that an increase on the values of $\tau_1$, $\tau_2$ and $\tau_3$ 
from $0.5$ to $1$ day implies an increase in the amplitude of all 
the populations with a translation.}
% -----------------------------------
\vspace{-8pt}

% ----------------------------------

\section{Conclusions and Discussion}
\label{sec:7}

In this paper, we have proposed and analyzed a three-compartment model 
that depicts the interaction between two prey and one predator. 
The new stochastic predator–prey model incorporates L\'{e}vy noise 
and considers time delays of the two prey in the logistic functional production 
and in the transformation rate of predator to prey. Both white noise 
and L\'{e}vy jump perturbations are integrated into all model compartments. We have 
established the existence and uniqueness of a global positive solution and its boundedness, 
demonstrating the well-posedness of our stochastic predator–prey mathematical model. 
We have also provided sufficient conditions for the extinction of both prey and predator 
as well as for the persistence of the prey with the extinction of the predator. 
Additionally, we have presented a sufficient condition for the persistence in mean 
of both prey and the predator. Our theoretical findings are reinforced by different 
numerical results, demonstrating three possible scenarios: 
(i) the extinction of both prey and predator populations, 
(ii) the persistence of the prey and the extinction of the predator, and 
(iii) biological persistence, which is vital for maintaining the continuity 
of the reaction between prey and predator.

{Several interesting research questions require further investigation in future studies. 
Specifically, our team intends to examine the spatial diffusion of a delayed prey–predator problem 
with L\'{e}vy jumps. Additionally, we plan to explore the bifurcation 
and chaotic behavior of the stochastic time-delayed prey–predator model. A comprehensive analysis 
is nontrivial and requires further inquiry, which we will undertake in another study. 
Furthermore, our analysis of the model under consideration can be expanded to 
include fractional-order derivative models, as described in \cite{34,35}, and we can explore 
alternative views of stochasticity, as demonstrated by the authors of \cite{36}. 
We also note that the parameter values presented in Table~\ref{tabl1} for numerical simulations 
of our theoretical results may not be applicable to real-world problems. 
Future research should strive to establish a connection with practical objectives 
and address the challenge of parameter estimation.}

% ----------------------------------

\vspace{6pt} 

\authorcontributions{Conceptualization, J.D. and D.F.M.T.; 
methodology, J.D. and D.F.M.T.; 
software, J.D.; 
validation, J.D. and D.F.M.T.; 
formal analysis, J.D. and D.F.M.T.; 
investigation, J.D. and D.F.M.T.; 
writing---original draft preparation, J.D. and D.F.M.T.; 
writing---review and editing, J.D. and D.F.M.T.; 
visualization, J.D.; 
project administration, D.F.M.T. 
All authors have read and agreed 
to the published version of the manuscript.}

\funding{This research was funded by 
The Portuguese Foundation for Science and Technology (FCT) 
grant numbers UIDB/04106/2020 and UIDP/04106/2020.}

\institutionalreview{Not applicable.}

\informedconsent{Not applicable.}
	
\dataavailability{No new data were created or analyzed in this study. 
Data sharing is not applicable to this article.} 

\acknowledgments{The authors are very grateful 
to three anonymous referees for reading carefully the submitted
manuscript and providing several valuable suggestions and comments.}

\conflictsofinterest{The authors declare no conflict of interest.
The funders had no role in the design of the study; in the collection, 
analyses, or interpretation of data; in the writing of the manuscript; 
or in the decision to publish the~results.} 

% ----------------------------------

\begin{adjustwidth}{-\extralength}{0cm}

\reftitle{References}
	
% ----------------------------------

% ----------------------------------

\PublishersNote{}

\end{adjustwidth}

\end{document}